 \documentclass[12pt]{article}
\usepackage{amsmath,amsthm,amssymb}

\newtheorem{thm}{Theorem}
\newtheorem{prop}{Proposition}
\newtheorem{lemma}{Lemma}

\newcommand{\ga}{\alpha}
\newcommand{\gb}{\beta}

\newcommand{\gd}{\delta}
\newcommand{\gep}{\epsilon}

\newcommand{\gs}{\sigma}

\newcommand{\gra}{\nabla}
\newcommand{\de}{\partial}
\newcommand{\bpf}{\begin{proof}}
\newcommand{\epf}{\end{proof}}
\newcommand{\beq}{\begin{equation}}
\newcommand{\eeq}{\end{equation}}
\newcommand{\beqn}{\begin{eqnarray*}}
\newcommand{\eeqn}{\end{eqnarray*}}

\begin{document}
\title{Conformal Deformation to Scalar Flat Metrics with Constant Mean Curvature on the Boundary in Higher Dimensions} \vskip 1em
 \author{ Szu-yu Sophie Chen
 \footnote{
 The author was supported by NSF grant DMS-0635607.}}
\date{}
\maketitle


  On a closed Riemannian manifold of dimension $n \geq 3$, every metric is conformal to  
  a constant scalar curvature metric. 
   This problem, called the Yamabe problem, was proved by Yamabe \cite{Yam60},
 Trudinger \cite{Tr68}, Aubin \cite{Aub76} and Schoen \cite{Sch84}.

 To extend the conformal deformation problem to manifolds with boundary,
  Escobar proposed two types of formulations.
 Let $(M, g)$ be a compact Riemannian manifold of dimension
 $n \geq 3$ with boundary $\de M.$
  We denote by $R_g$ the scalar curvature of the manifold
   and by $\kappa_g$ the mean curvature of the boundary.
 The first type is to 
 find a metric $\tilde g$ in the conformal class of $g$
    such that $R_{\tilde g}$ is constant and $\kappa_{\tilde g}$ is
    zero.  This  was  studied by Escobar \cite{Es92} and recently
  by Brendle and the author \cite{Bchen09}. 
  
   The second type is to find a metric $\tilde g$ in the conformal class of $g$
    such that $R_{\tilde g}$ is zero and $\kappa_{\tilde g}$ is
    constant. This problem, as Escobar remarked  \cite{Es92a}, is  a higher
  dimensional generalization of the Riemann mapping theorem.  The problem
  is studied  by Escobar \cite{Es92a}, \cite{Es96} and Marques \cite{Marques05}, \cite{Marques07}.
   (For analysis background for both problems, see \cite{Cherrier84}).

  In this paper, we will study the second formulation; that is the  existence of a conformal metric with zero scalar curvature and constant mean 
  curvature on the boundary. The
  problem turns out to be finding a critical point of the
  functional
  $$E_g(\phi) = \frac{\int_M (\frac{4(n-1)}{n-2} |\gra_g \phi|^2 + R_g \phi^2) dV_g
  + \int_{\de M}  2  \kappa_g  \phi^2 d\gs_g}{(\int_{\de M} \phi^{\frac{2(n-1)}{n-2}}d\gs_g)^{\frac{n-2}{n-1}}},$$
  where $\phi$ is a  positive smooth function on $M.$
  The exponent $\frac{2(n-1)}{n-2}$ is critical for  the trace Sobolev
  embedding $H^1 (M) \hookrightarrow L^{\frac{2(n-1)}{n-2}} (\de M).$
  This embedding is not compact and the functional $E_g$
   does not satisfy the Palais-Smale condition. For this reason, standard
   variational methods cannot be applied.

   To study the problem, we consider
   the Sobolev quotient, introduced in \cite{Es92a},
  $$\mathcal{Q}(M, \de M, g) = \inf_{0< \phi \in C^{\infty}} E_g(\phi).$$
  This is known that $\mathcal{Q} (M, \de M, g)$ is a conformal invariant and
   $\mathcal{Q}(M, \de M, g) \leq \mathcal{Q}(B^n, \de B^n),$
  where $\mathcal{Q}(B^n, \de B^n)$ is the Sobolev quotient of the
  unit ball $B^n$ in $\mathbb{R}^n$ equipped with the flat metric.
  It  was proved by Escobar that

 \begin{thm} (Escobar \cite{Es92a}) If $\mathcal{Q}(M, \de M, g) < \mathcal{Q}(B^n, \de B^n),$
 then there exists a metric $\tilde g$ in the conformal class of $g$
 such that $R_{\tilde g}$ is zero and $\kappa_{\tilde g}$ is
 constant.
 \end{thm}

 For $n \geq 6,$ when $\de M$ is not umbilic,
 Escobar showed that $\mathcal{Q}(M, \de M, g) < \mathcal{Q}(B^n, \de B^n).$
   He also proved the inequality holds when $n=3,$ and when $n=4, 5$ and
  $\de M$ is umbilic, provided $M$ is not
  conformally equivalent to the unit ball.
  When  $n=4, 5,$ and  $\de M$ is not umbilic, Marques
  verified that the inequality holds.

 Consequently, it remains to consider the case that $n \geq 6$ and $\de M$ is
 umbilic (some special case was considered in \cite{Marques05}).
 As in \cite{Brendle07}, \cite{Bchen09}, we denote by
 $\mathcal{Z}$ the set of points $p \in M$ such that
    $$\limsup_{x \rightarrow p} d(p, x)^{2-d} |W_g| (x) = 0,$$
  where $d = [\frac{n-2}{2}]$ and $W_g$ is the Weyl tensor of $g.$
  We note that $p \in \mathcal{Z}$ if  and only if $\gra^m W_g (p)= 0$
  for $m = 0, \cdots, d-2.$ Moreover, the set $\mathcal{Z}$ is
  conformally invariant.

 Our main result  is
 \begin{thm} \label{t:main} Let $(M, g)$ be a compact Riemannian
 manifold of dimension $n \geq 6$ with umbilic boundary. Suppose there
 exists a point $p \in \de M$ such that $p \notin \mathcal{Z},$ then
   $\mathcal{Q}(M, \de M, g) < \mathcal{Q}(B^n, \de B^n).$ As a
   result, there exists a metric $\tilde g$ in the conformal class of $g$
 such that $R_{\tilde g}$ is zero and $\kappa_{\tilde g}$ is
 constant.
 \end{thm}

We now discuss the case that $p \in \mathcal{Z}$ for all $p \in \de M.$
 In Section~\ref{s:proof}, we consider a flux integral $\mathcal{I} (p, \gd)$
  introduced in \cite{Bchen09}
 in a small neighborhood of $p\in \de M.$ When $p \in \mathcal{Z},$
 it was shown in \cite{Bchen09} that $\lim_{\gd \rightarrow 0} \mathcal{I} (p, \gd)$
 exists and is equal to a positive multiple of ADM mass of  certain scalar flat
 asymptotically flat manifold; see  Section~\ref{s:proof}.
 We reduce the case to positivity of mass.

 \begin{thm} \label{t:main2} Let $(M, g)$ be a compact Riemannian
 manifold of dimension $n \geq 6$ with umbilic boundary. Suppose there
 exists a point $p \in \de M$ such that $p \in \mathcal{Z}$ and
 $\lim_{\gd \rightarrow 0} \mathcal{I} (p, \gd)> 0,$ then
   $\mathcal{Q}(M, \de M, g) < \mathcal{Q}(B^n, \de B^n).$
As a result, there exists a metric $\tilde g$ in the conformal class
of $g$
 such that $R_{\tilde g}$ is zero and $\kappa_{\tilde g}$ is
 constant.
 \end{thm}

 We give the outline of the proof. By Marques \cite{Marques05}, we may
 choose  conformal Fermi coordinates around a boundary point $p.$
 In these coordinates, we define
   $$v_{\gep} =  \left(\frac{\gep}{(\gep+ x_n)^2 +\sum_{1 \leq a \leq n-1} x_a^2}\right)^{
   \frac{n-2}{2}}.$$
  We note that $v_{\gep}$ is the extremal function for the sharp trace
  Sobolev inequality on the half plane; see \cite{Es88}, \cite{Beckner93}.
  By conformal invariance, it holds
   $$\mathcal{Q}(B^n, \de B^n) \left( \int_{\de \mathbb R^n_+}  v_{\gep}^{\frac{2(n-1)}{n-2}} d \gs\right)^{\frac{n-2}{n-1}} = \frac{4(n-1)}{n-2}\int_{\mathbb R^n_+} |\gra v_{\gep}|^2 dx.$$
  It is then understood that  $v_{\gep}$ is the model function on $\mathbb R^n_+.$
  
   We now consider the function $v_{\gep} + \psi$ defined in a small neighborhood of
   $p,$  where $\psi$ satisfies
    \begin{align}
  \Delta \psi =  \sum_{i,k = 1}^n (\frac{n-2}{4(n-1)} v_{\epsilon} \de_i \de_k S_{ik} + \de_k (\de_i v_{\epsilon} S_{ik})) & \quad \text{in} \; B_{\gd} \cap \mathbb{R}^n_+, \label{e:eqn} \\
  \de_n \psi = -\frac{1}{2 (n-1)} \de_n v_{\gep} S_{nn} + \frac{n}{n-2}
  v_{\gep}^{-1}
  \de_n v_{\gep} \psi &\quad \text{on}\; B_{\gd} \cap \de \mathbb{R}^n_+.
  \label{e:bdy}
  \end{align} 
  In the above equations,
  the tensor $S_{ij}$ comes from applying  the conformal killing operator to certain vector
  field we solve; see   Section~\ref{s:construct}.
  The equation (\ref{e:eqn})
  corresponds to a  linear  approximation of the scalar curvature equation of $(v_{\gep} + \psi)^{\frac{4}{n-2}} g.$
  However, in our construction, the boundary condition (\ref{e:bdy})
   is  not the
  linear approximation of the mean curvature equation of $(v_{\gep} + \psi)^{\frac{4}{n-2}} g;$
  the "linear mean curvature equation" should be
   $$\de_n \psi =  \frac{n}{n-2}  v_{\gep}^{-1}  \de_n v_{\gep} \psi.$$ 
  We emphasize that the Sobolev quotient $\mathcal{Q} (M, \de M, g)$
  is  normalized  by the volume of the boundary (not the volume of the manifold).
   Our deformation of the metric does not fix the volume of the boundary locally.
   As a consequence, in order to get the energy functional small enough, 
   the term $-\frac{1}{2 (n-1)} \de_n v_{\gep} S_{nn}$ is important because it 
     cancels out \emph{to the right order}  the change of  the volume of the boundary. 
   This is the reason that the linear approximation of the mean curvature equation does not work here.
    This turns out to be the delicate part of the proof.
   Finally, to define a test function globally, we glue the
   function $v_{\gep} + \psi$ with the Green's function of the
   conformal Laplacian  centered at $p.$

To show  the above test function has the energy functional less
than $\mathcal{Q} (B^n , \de B^n),$ we use  the method and
techniques developed by Brendle \cite{Brendle07} (see also
\cite{Bchen09}). In \cite{Brendle07}, these nice techniques were
used
 to prove a convergence theorem for the Yamabe flow. 
  In \cite{Bchen09}, these techniques were used to study the problem of  first type 
 described at the beginning.
 To be more precise, let $u_{\gep} = \gep^{\frac{n-2}{2}} (\gep^2 + |x|^2)^{-\frac{n-2}{2}}.$
 In \cite{Brendle07}, one considers the function $u_{\gep} + w$ in
 normal coordinates,
 where $w$ satisfies
  $
  \Delta w + n(n+2) u_{\epsilon}^{\frac{4}{n-2}} w = \frac{n-2}{4(n-1)} u_{\epsilon} \de_i \de_k S_{ik} + \de_k (\de_i u_{\epsilon} S_{ik}).$
  In \cite{Bchen09}, one considers the function $u_{\gep} + w$ in
 Fermi coordinates together with the boundary condition $\de_n w = 0.$
 We refer the readers to \cite{Brendle05}, \cite{Brendle08},
 \cite{KMS07}, \cite{BM09}, \cite{Chen05b} for other related works concerning the Yamabe
 problem.

   We introduce the  notation in this paper. We denote by $dx$ the volume element in $\mathbb{R}^n,$
  by $d\gs$ the area element of a hypersurface in $\mathbb{R}^n$ and
  by $d \mu$ the area element of an $(n-2)$-dimensional surface in
  $\mathbb{R}^n.$ We also denote by $\mathbb{R}^n_+$ the half plane
  $\{x: x_n \geq 0\}.$
  Let $B_r(x)$ be the ball of radius $r$ centered at $x.$ When $x$
  is at the origin, we simply denote by $B_r.$

  \textbf{Acknowledgment:}  The author would like to thank Simon Brendle for valuable comments which help
   improve the presentation of the work.

\section{Background} \label{s:variation}

 Let
  $v_{\gep} (x) = \gep^{\frac{n-2}{2}} ((\gep+ x_n)^2 +\sum_{a= 1}^{n-1}x_a^2)^{- \frac{n-2}{2}}$,
   $x \in \mathbb{R}^n_+.$
  The
  function  $v_{\gep}$ satisfies
  \beq \label{e:v}
 \Delta v_{\gep} = 0 \quad\text{for}\; x \in \mathbb{R}^n_+,
  \eeq
  \beq \label{e:hess_v}
  v_{\gep} \de_i \de_k v_{\gep} - \frac{n}{n-2} \de_i v_{\gep} \de_k v_{\gep}=
   - \frac{1}{n-2} |d v_{\gep}|^2 \gd_{ik} \quad \text{for}\; x \in \mathbb{R}^n_+,
  \eeq  and
  \beq \label{e:v_bdy}
  \de_n v_{\gep}= - (n-2) v_{\gep}^{\frac{n}{n-2}}
  \quad \text{for}\; x \; \in \de \mathbb{R}^n_+.
  \eeq
  By integration, we get
  $$\int_{\mathbb R^n_+} |\gra v_{\gep}|^2 dx = (n-2) \int_{\de \mathbb R^n_+} v_{\gep}^{\frac{2(n-1)}{n-2}} d\gs.$$
  Moreover,   $v_{\gep}$ satisfies the following
  inequalities:
   $$\gep^{\frac{n-2}{2}} (\gep + |x|)^{-n+2} \leq v_{\gep} (x) \leq C(n) \gep^{\frac{n-2}{2}} (\gep + |x|)^{-n+2}
   \quad \text{for}\; x \in \mathbb{R}^n_+;$$
   $$|\de v_{\gep}| (x) \leq C(n) \gep^{\frac{n-2}{2}} (\gep +
   |x|)^{-n+1} \quad \text{for}\; x \in \mathbb{R}^n_+;$$
  and 
   $$|v_{\gep} - \gep^{\frac{n-2}{2}} |x|^{-n+2}| \leq C(n) \gep^{\frac{n}{2}} |x|^{-n+1}
    \quad \text{for}\;  x \in \mathbb{R}^n_+,\; \text{and}\; |x|  \geq 2\gep,$$
   where $C(n)$ is a positive constant depending only on $n.$

  Let $V$ be a smooth vector field and $H_{ik}$ be a trace-free symmetric two-tensor. We define
     $$\begin{array} {l}
  S_{ik}= \de_i V_k + \de_k V_i - \frac{2}{n} div V \gd_{ik},\\
   T_{ik} = H_{ik}- S_{ik}, \\
    P_{ik,l} = v_{\gep} \de_l T_{ik} - \frac{2}{n-2} \de_i v_{\gep} T_{kl} - \frac{2}{n-2} \de_k v_{\gep} T_{il} + \frac{2}{n-2} \sum_{p=1}^n \de_p v_{\gep} T_{ip} \gd_{kl}
    + \frac{2}{n-2} \sum_{p=1}^n  \de_p v_{\gep} T_{kp} \gd_{il},\\
   \psi = \de_l v_{\epsilon} V_l + \frac{n-2}{2n} v_{\epsilon} div V.
   \end{array}$$
   In \cite{Brendle07}, \cite{Bchen09}, a similar notation was introduced with $v_{\gep}$  replaced by 
     $u_{\gep}.$

   The following formula    is a revision of the formula in \cite{Brendle07} Proposition~5, 6.
   The formula in \cite{Brendle07}  corresponds to   the second
   variation of the scalar curvature  on the sphere.
   Similarly, the formula here corresponds to  the  second variation of the scalar curvature on the ball in $\mathbb R^n$.
  \begin{prop}  \label{p:b1} Let $H_{ik}$ be a trace-free symmetric two-tensor, and  $V$ be a smooth vector field.
  Then $\psi$ satisfies
  \beq \label{e:Delta_w}
  \Delta \psi = \sum_{i,k = 1}^n (\frac{n-2}{4(n-1)} v_{\epsilon} \de_i \de_k S_{ik} + \de_k (\de_i v_{\epsilon} S_{ik})).
  \eeq
  Moreover,
  \begin{align}
   & \frac{1}{4} |P|^2 - \frac{1}{2} \sum_{i= 1}^n|\sum_{k= 1}^n (v_{\gep} \de_k T_{ik}+ \frac{2n}{n-2} \de_k v_{\gep} T_{ik})|^2
       \notag\\
    =&\, \sum_{i,k, l = 1}^n (\frac{1}{4} v_{\gep}^2 \de_l H_{ik} \de_l H_{ik} - \frac{1}{2} v_{\gep}^2 \de_k H_{ik} \de_l H_{il} -2 v_{\gep} \de_k v_{\gep} H_{ik} \de_l H_{il}
           - \frac{2(n-1)}{n-2}  \de_k v_{\epsilon} \de_l v_{\epsilon} H_{ik} H_{il}) \notag\\
      + & \sum_{i,k = 1}^n (- 2 v_{\gep} \psi \de_i \de_k H_{ik} +  \frac{8 (n-1)}{n-2} \de_i v_{\gep} \de_k \psi  H_{ik}) - \frac{4(n-1)}{n-2} |d\psi|^2
        + \sum_{i=1}^n \de_i \xi_i , \notag
     \end{align}
     where
       \begin{align}
       \xi_i
       = & \, \sum_{k= 1}^n (2 v_{\gep} \psi \de_k H_{ik} - 2v_{\gep} \de_k \psi H_{ik} - 2\de_k v_{\gep} \psi H_{ik}
         -v_{\gep} \psi \de_k S_{ik} + \de_k (v_{\gep} \psi) S_{ik})  \notag\\
      + & \sum_{k, l = 1}^n  (2 v_{\gep} \de_l v_{\gep} S_{kl} H_{ki} - \frac{1}{2} v_{\gep}^2 \de_i S_{lk} H_{lk}
      +  v_{\gep}^2 \de_l S_{kl} H_{ki} +\frac{1}{4}  v_{\gep}^2 \de_i S_{lk} S_{lk}
       - \frac{1}{2}  v_{\gep}^2 \de_k S_{lk} S_{il})\notag \\
      + &  \sum_{k, l =1}^n (-  v_{\gep} \de_k  v_{\gep} S_{lk} S_{il}
        - \frac{2}{n-2} v_{\gep} \de_k v_{\gep} T_{lk} T_{li})
        + \frac{4(n-1)}{n-2} (- \sum_{k= 1}^n\de_k v_{\gep} \psi S_{ik} + \psi \de_i
        \psi). \notag
        \end{align}
 \end{prop}
  \bpf Since the proof is similar, we only point out the difference.
   In \cite{Brendle07} Proposition 5, it was shown that
    \beqn v_{\epsilon} \de_i \de_k S_{ik} + \frac{4 (n-1)}{n-2} \de_k (\de_i v_{\epsilon} S_{ik})
    &=& \frac{4 (n-1)}{n-2} \Delta (\sum_{l=1}^n\de_l v_{\epsilon} V_l + \frac{n-2}{2n} v_{\epsilon} div
    V)\\
    &-& \frac{4 (n-1)}{n-2} (\sum_{l= 1}^n\de_l \Delta v_{\epsilon} V_l + \frac{n+2}{2n} \Delta v_{\epsilon} div V)
    \eeqn (with $v_{\gep}$ replaced by $u_{\gep}$ but the formula holds in
    general).
    By (\ref{e:v}), then (\ref{e:Delta_w}) follows.

    For the second identity,
    by \cite{Brendle07} Proposition 5, it holds
    \beqn
    & &\frac{1}{4} v_{\epsilon}^2 |\de T|^2 - \frac{1}{2} v_{\gep}^2 |div T|^2 - \sum_{i,k,l= 1}^n (2
    v_{\epsilon} \de_k v_{\epsilon} T_{ik} \de_l T_{il} +
    \frac{2(n-1)}{n-2} \de_k v_{\epsilon} \de_l v_{\epsilon} T_{ik}
    T_{il})\\
    &=& I_1 - 2 I_2 + I_3,
    \eeqn where
    $$ 
    I_1 =  \sum_{i,k, l = 1}^n (\frac{1}{4} v_{\gep}^2 \de_l H_{ik} \de_l H_{ik} - \frac{1}{2} v_{\gep}^2 \de_k H_{ik} \de_l H_{il} -2 v_{\gep} \de_k v_{\gep} H_{ik} \de_l H_{il}
           - \frac{2(n-1)}{n-2}  \de_k v_{\epsilon} \de_l v_{\epsilon} H_{ik}
           H_{il}),
    $$
    \beqn
     I_2 &=& \sum_{i, k = 1}^ n (v_{\gep} \psi \de_i \de_k H_{ik}- \frac{4(n-1)}{n-2} \de_i
     v_{\gep} \de_k \psi H_{ik} - \de_i ( v_{\gep} \psi \de_k H_{ik}) +
     \de_k (v_{\gep} \de_i \psi H_{ik}))\\
     & +& \sum_{i,k = 1}^n \de_k (\de_i v_{\gep} \psi H_{ik}) + \sum_{i,k,l =1}^n (\frac{1}{4} \de_l (v_{\gep}^2 \de_l S_{ik}
     H_{ik})-\frac{1}{2} \de_k (v_{\gep}^2 \de_l S_{il} H_{ik})
     -\de_k (v_{\gep} \de_l v_{\gep} S_{il} H_{ik}))\\
     & +& \sum_{i,k,l =1}^n (v_{\gep} \de_k \de_l v_{\gep} - \frac{n}{n-2} \de_k v_{\gep} \de_l
     v_{\gep}) (\de_l V_i - \de_i V_l) H_{ik}\\
     &-& \sum_{i,k,l = 1}^n \de_l [(v_{\gep} \de_i \de_k v_{\gep} - \frac{n}{n-2} \de_i v_{\gep} \de_k
     v_{\gep}) V_l] H_{ik},
    \eeqn and
    \beqn
     I_3 &=& \sum_{i, k = 1}^ n (v_{\gep} \psi \de_i \de_k S_{ik}- \frac{4(n-1)}{n-2} \de_i
     v_{\gep} \de_k \psi S_{ik} - \de_i ( v_{\gep} \psi \de_k S_{ik}) +
     \de_k (v_{\gep} \de_i \psi S_{ik}))\\
     & +& \sum_{i,k = 1}^n \de_k (\de_i v_{\gep} \psi S_{ik}) + \sum_{i,k,l =1}^n (\frac{1}{4} \de_l (v_{\gep}^2 \de_l S_{ik}
     S_{ik})-\frac{1}{2} \de_k (v_{\gep}^2 \de_l S_{il} S_{ik})
     -\de_k (v_{\gep} \de_l v_{\gep} S_{il} S_{ik}))\\
     & +& \sum_{i,k,l =1}^n (v_{\gep} \de_k \de_l v_{\gep} - \frac{n}{n-2} \de_k v_{\gep} \de_l
     v_{\gep}) (\de_l V_i - \de_i V_l) S_{ik}\\
     &-& \sum_{i,k,l = 1}^n \de_l [(v_{\gep} \de_i \de_k v_{\gep} - \frac{n}{n-2} \de_i v_{\gep} \de_k
     v_{\gep}) V_l] S_{ik}.
    \eeqn  And in \cite{Brendle07} Proposition 6, it holds
  \begin{align}
   &  \frac{1}{4} v_{\epsilon}^2 |\de T|^2 - \frac{1}{2} v_{\gep}^2 |div T|^2 - \sum_{i,k,l= 1}^n (2
    v_{\epsilon} \de_k v_{\epsilon} T_{ik} \de_l T_{il} +
    \frac{2(n-1)}{n-2} \de_k v_{\epsilon} \de_l v_{\epsilon} T_{ik}
    T_{il}) \notag\\
    &=  \frac{1}{4} |P|^2
    - \frac{1}{2} \sum_{i= 1}^n|\sum_{k= 1}^n (v_{\gep} \de_k T_{ik}
    + \frac{2n}{n-2} \de_k v_{\gep} T_{ik})|^2- \frac{2}{(n-2)^2} |\de
    v_{\gep}|^2 |T|^2  \notag\\
    &+ \sum_{i,k,l= 1}^n  (- \frac{2}{n-2} (v_{\gep} \de_k \de_l v_{\gep} - \frac{n}{n-2} \de_k v_{\gep} \de_l
    v_{\gep}) T_{ik} T_{il} + \frac{2}{n-2} \de_l (v_{\gep} \de_k v_{\gep} T_{ik}
    T_{il})) \label{e:I}
  \end{align} (with $v_{\gep}$ replaced by $u_{\gep}$ but the formula holds in
  general).
 Using  (\ref{e:hess_v}) in $I_2$,  (\ref{e:hess_v}) and (\ref{e:Delta_w}) in $I_3$ and
     using (\ref{e:hess_v}) in (\ref{e:I}) give the identity.     
  \epf

\section{Construction} \label{s:construct}
   We first state some properties about conformal
    Fermi coordinates that we will use later.
   Then we construct the correction term $\psi$ and compute some formulas
  on the boundary. Let $n \geq 6.$ We assume $\de M$ is totally geodesic.

  In this section, we assume $g$ is the metric in conformal Fermi
  coordinates. We write  $g = \exp h.$
 By Marques \cite{Marques05}, 
   we have  $tr\, h(x) = O(|x|^{2d+2})$
 for $x \in \mathbb{R}^n_+,$ where $d = [\frac{n-2}{2}].$
 Moreover,  
    $h_{i n}(x) = 0$  for $x \in \mathbb{R}^n_+$ and $i = 1, \cdots, n.$  
We also have   $\de_n h_{ab} (x)= \sum_{i=1}^n h_{ai}(x) x_i = 0$  for $x \in \de \mathbb{R}^n_+$ 
 and $a, b = 1, \cdots, n-1.$  
  In this case, $\det g (x)= 1+ O(|x|^{2d+2})$ for $x \in \mathbb{R}^n_+.$
 

  Let  $H_{ij}$ be the Taylor expansion of  $h_{ij}$ up to the order $d$
  $$
  H_{ij} = \sum_{2 \leq |\ga| \leq d} h_{ij, \ga} (0) x^{\ga},
  $$ where $\ga$ is a multi-index. Then $h_{ik} = H_{ik} + O(|x|^{d+1}).$
  It follows that
  $$tr H (x)= H_{i n} (x)= 0 \quad  \text{for all}  \; x \in \mathbb R^n_+ \;\text{and}\; i= 1, \cdots, n,$$  and
   $$\de_n H_{ab}(x)= \sum_{i= 1}^n H_{ai}(x) x_i = 0  \quad \text{for all} \; x \in \de \mathbb{R}^n_+
   \;\text{and}\;
   a, b= 1, \cdots, n-1.$$

   We define algebraic Schouten tensor and algebraic Weyl tensor of $H_{ij}$ as in \cite{Brendle07}:
  $$\begin{array} {l}
   A_{ij}= \de_i \de_m H_{mj} + \de_m \de_j H_{im} - \Delta H_{ij} - \frac{1}{n-1} \de_m \de_p H_{mp} \gd_{ij},\\
   Z_{ijkl} = \de_i \de_k H_{jl} - \de_i \de_l H_{jk} -\de_j \de_k H_{il} + \de_j \de_l H_{ik} + \frac{1}{n-2} (A_{jl} \gd_{ik} - A_{jk} \gd_{il}- A_{il} \gd_{jk}+ A_{ik} \gd_{jl}).
  \end{array}$$

    \begin{prop} \label{p:fermi} \cite{Bchen09}   If  $Z_{ijkl}= 0$ for all $x \in \mathbb R^n_+,$
  then $H_{ij} = 0$ for all $x \in \mathbb R^n_+.$
  \end{prop}


 \begin{prop} \label{p:scalar} \cite{Bchen09} The scalar curvature
 $R_g$ satisfies
  $$|R_g - \de_i \de_k H_{ik}| \leq C \sum_{i,j}  \sum_{2 \leq |\ga| \leq d} |h_{ij, \ga}| |x|^{|\ga|} + C |x|^{d-1},$$ and
  \beqn
  & &|R_g - \de_i \de_k h_{ik} + \de_k (H_{ik} \de_l H_{il}) - \frac{1}{2} \de_k H_{ik} \de_l H_{il} + \frac{1}{4} \de_l H_{ik} \de_l H_{ik}|\\
  &\leq& C \sum_{i, j}  \sum_{2 \leq |\ga| \leq d} |h_{ij, \ga}|^2 |x|^{2 |\ga|}  + C \sum_{i, j}  \sum_{2 \leq |\ga| \leq d} |h_{ij, \ga}| |x|^{|\ga| +d  -1} + C |x|^{2 d}
  \eeqn
  for $|x|$ sufficiently small.
  \end{prop}

    Let $V$ be a smooth vector field.
    We next define as in Section~\ref{s:variation} that
     $$\begin{array} {l}
  S_{ik}= \de_i V_k + \de_k V_i - \frac{2}{n} div V \gd_{ik},\\
   T_{ik} = H_{ik}- S_{ik}, \\
    P_{ik,l} = v_{\gep} \de_l T_{ik} - \frac{2}{n-2} \de_i v_{\gep} T_{kl} - \frac{2}{n-2} \de_k v_{\gep} T_{il} + \frac{2}{n-2} \sum_{p=1}^n  \de_p v_{\gep} T_{ip} \gd_{kl}
    + \frac{2}{n-2}\sum_{p=1}^n \de_p v_{\gep} T_{kp} \gd_{il}.
   \end{array}$$
    \begin{prop} \label{p:P} Let $V$ be a smooth vector field. Then
    $$\sum_{i, j}  \sum_{2 \leq |\ga| \leq d} |h_{ij, \ga}|^2 \gep^{n-2} \int_{B_{\gd} \cap \mathbb{R}^n_+} (\gep + |x|) ^{2 |\ga| +2 -2n} dx  \leq C(n)
      \int_{B_{\gd} \cap \mathbb{R}^n_+} |P|^2 dx$$
    for $\gd \geq 2 \gep > 0.$
 \end{prop}
  \bpf In \cite{Brendle07} Proposition 9, it was shown that
  \begin{align}
  &\sum_{i,j,k,l =1}^n \{\de_j (\de_l T_{ik} - \frac{2}{n-2}v_{\gep}^{-1} \de_k v_{\gep} T_{il})
   + \frac{2}{n-2} v_{\gep}^{-1} \de_k v_{\gep}
   (\de_j T_{il} -\frac{2}{n-2} v_{\gep}^{-1} \de_i v_{\gep} T_{jl})
   \notag\\
 & + \frac{2}{n-2} v_{\gep}^{-2} (v_{\gep} \de_j \de_k v_{\gep} -\frac{n}{n-2} \de_j v_{\gep} \de_k v_{\gep})
   T_{il} + \frac{4}{(n-2)^2} v_{\gep}^{-2} \de_k v_{\gep} (\de_i v_{\gep} T_{jl} + \de_j v_{\gep}
   T_{il})\} Z_{ijkl}\notag\\
   & = \sum_{i,j,k,l =1}^n \de_j \de_l H_{ik} Z_{ijkl} \notag
  \end{align} (with $v_{\gep}$ replaced by $u_{\gep}$ but the formula holds in
  general).
  Then by (\ref{e:hess_v}), we have
   $$\sum_{i,j,k,l = 1}^n (\de_j (v_{\gep}^{-1} P_{ik,l}) Z_{ijkl} + \frac{2}{n-2} v_{\gep}^{-2} \de_k v_{\gep}
   P_{il,j} Z_{ijkl} )= \frac{1}{4}|Z|^2.$$
  From this, the assertion follows easily by the proof in \cite{Bchen09} Proposition 7 and Corollary 8
  using $\gep^{\frac{n-2}{2}} (\gep + |x|)^{-n+2} \leq v_{\gep} (x) \leq C(n) \gep^{\frac{n-2}{2}} (\gep + |x|)^{-n+2}$
  and  $|\de v_{\gep}| (x) \leq C(n) \gep^{\frac{n-2}{2}} (\gep +
   |x|)^{-n+1}.$
  \epf

  We next construct the correction term $\psi$.
  We fix a  positive smooth function $\eta (t)$ such that $\eta(t)
  =1$ for $t \leq \frac{4}{3}$ and $\eta(t) = 0$ for $t \geq \frac{5}{3}.$
  For $\gd > 0,$ we define $\eta_{\gd} (x) = \eta (\frac{|x|}{\gd})$,
   $x \in \mathbb{R}^n_+.$
  Notice that $\de_n \eta_{\gd} (x) = 0$ for all $x \in \de \mathbb{R}^n_+.$
  By Proposition~\ref{p:system} in Appendix,
   there exists  a smooth vector field $V$  which solves
     \beq \label{e:star} \left\{  \begin{array}{ll}
  \sum_{k =1}^n \de_k [ v_{\gep}^{\frac{2n}{n-2}} ( \eta_{\gd} H_{ik} - \de_i V_k - \de_k V_i + \frac{2}{n} div V \gd_{ik})]= 0 & in \, \mathbb{R}^n_+  \\
  \de_n V_a=  0 & on \, \de \mathbb{R}^n_+\\
  V_n = 0 & on \, \de \mathbb{R}^n_+
  \end{array}\right.
    \eeq for $i = 1, \cdots, n$ and $a = 1, \cdots, n-1.$
   Moreover, $V$ satisfies
    \beq \label{i:V}
    |\de^{\gb} V^{(\gep, \gd)}(x)|  \leq C(n, |\gb|) \sum_{i,k= 1}^n \sum_{2 \leq |\ga| \leq d} |h_{ik, \ga}| (\gep+ |x|)^{|\ga|+1 - |\gb|}.
    \eeq
   By the equation,
     \beq \label{e:T}
    \sum_{k= 1}^n (v_{\gep} \de_k T_{ik}+ \frac{2n}{n-2} \de_k v_{\gep} T_{ik}) = 0
    \eeq for $x \in B_{\gd} \cap \mathbb{R}^n_+$ and $i = 1, \cdots, n.$
    We next define
    $$\psi =\sum_{l=1}^n \de_l v_{\epsilon} V_l + \frac{n-2}{2n} v_{\epsilon} div
    V.$$

 \begin{prop} \label{p:bdy}  It holds
  $ S_{a n} (x)= 0$,
   $$\de_n S_{nn}(x)= - \frac{2n}{n-2} v_{\gep}(x)^{-1} \de_n
  v_{\gep} (x)  S_{nn} (x)=   2n\, v_{\gep}(x)^{\frac{2}{n-2}}   S_{nn} (x),$$
  and
   $$\de_n S_{a b} (x)= - \frac{2n}{n-1} v_{\gep}(x)^{\frac{2}{n-2}} S_{nn}(x) \gd_{ab}$$
   for $ x \in \de \mathbb{R}^n_+$ and $ a, b = 1, \cdots, n-1.$
     As a consequence,  for $ x \in \de \mathbb{R}^n_+,$
  $$\de_n \psi (x)= -\frac{1}{2 (n-1)} \de_n v_{\gep}(x) S_{nn}(x) + \frac{n}{n-2}
  v_{\gep}(x)^{-1}
  \de_n v_{\gep}(x) \psi(x).$$
  \end{prop}
  \bpf
   By assumptions, $V_n= \de_n V_a = 0$ for  $ x \in \de
   \mathbb{R}^n_+$ and $a= 1, \cdots, n-1.$
   Thus, $S_{na} = T_{na} = \de_n V_a - \de_a V_n = 0$ on $\de \mathbb{R}^n_+$ for  $a= 1, \cdots, n-1$
   and   
   \beq\label{e:DnDaVb}
    \de_n \de_a V_b=0  \quad \text{for}\; x \in \de \mathbb{R}^n_+  \,\; \text{and}\; a, b= 1, \cdots, n-1.
     \eeq
     
   We next consider the equation (\ref{e:star}). It gives  $\sum_{k= 1}^n (v_{\gep} \de_k (\eta_{\gd}H_{nk}- S_{nk}) + \frac{2n}{n-2} \de_k v_{\gep} (\eta_{\gd} H_{nk}- S_{nk}))= 0.$ Since $H_{nk} (x) = 0$ for all  $x \in \mathbb{R}^n_+$ and  $k = 1, \cdots, n,$ we have
     $$\sum_{k= 1}^n (v_{\gep} \de_k  S_{nk} + \frac{2n}{n-2} \de_k v_{\gep}  S_{nk})= 0.$$
    for all  $x \in \mathbb{R}^n_+.$ Therefore, using (\ref{e:v_bdy})
   $$\de_n S_{nn} =  - \sum_{a=1}^{n-1} \de_a S_{na}- \frac{2n}{n-2} v_{\gep}^{-1} \sum_{k=1}^n \de_k v_{\gep}  S_{nk}
   = - \frac{2n}{n-2} v_{\gep}^{-1} \de_n v_{\gep} S_{nn} =   2n v_{\gep}^{\frac{2}{n-2}}   S_{nn} .$$
   Moreover, by (\ref{e:DnDaVb}), it follows that 
   \begin{align}
   \de_n S_{ab} &= \de_n \de_a V_b + \de_n \de_b V_a -\frac{2}{n} \de_n div V \gd_{ab}= 
   -\frac{2}{n} \de_n \de_n V_n  \gd_{ab}\notag\\
   & = - \frac{1}{n-1} (2 \de_n V_n - \frac{2}{n} \de_n div V) \gd_{ab}
   = - \frac{1}{n-1} \de_n S_{nn} \gd_{ab} = - \frac{2n}{n-1} v_{\gep}^{\frac{2}{n-2}} S_{nn}
   \gd_{ab}. \notag
   \end{align}

  We now compute $\de_n \psi.$  
   \begin{align}
     \de_n \psi    
   =&   \sum_{i=1}^n (\de_n \de_i v_{\gep} V_i + \de_i v_{\gep} \de_n
    V_i)
    + \frac{n-2}{2n} \de_n v_{\gep} div V + \frac{n-2}{2n} v_{\gep}
    \de_n div V \notag\\
   =& \sum_{i= 1}^n (\de_n \de_i v_{\gep} - \frac{n}{n-2} v_{\gep}^{-1} \de_i v_{\gep} \de_n v_{\gep}) V_i + \sum_{i=1}^n \de_i v_{\gep} \de_n V_i\notag\\
      + & \frac{n}{n-2} (\sum_{i=1}^n \de_i v_{\gep} V_i + \frac{n-2}{2n} v_{\gep} div
    V) v_{\gep}^{-1} \de_n v_{\gep}  - \frac{1}{n} div
    V \de_n v_{\gep} + \frac{n-2}{2n} v_{\gep} \de_n div V.\notag 
        \end{align}
      By (\ref{e:hess_v}) and $\de_n div V = \frac{n}{2(n-1)} \de_n S_{nn},$ we get
     \begin{align}
      \de_n \psi =& \frac{n}{n-2} (\sum_{i=1}^n \de_i v_{\gep} V_i + \frac{n-2}{2n} v_{\gep} div
    V) v_{\gep}^{-1} \de_n v_{\gep}  - \frac{1}{n} div
    V \de_n v_{\gep} \notag \\
     - & \frac{1}{n-2} |d v_{\gep}|^2 V_n + \sum_{i=1}^n \de_i v_{\gep} \de_n V_i + \frac{n-2}{4(n-1)} v_{\gep} \de_n S_{nn}. \notag
      \end{align}             
   Since  $\de_n S_{nn} = - \frac{2n}{n-2} v_{\gep}^{-1} \de_n v_{\gep} S_{nn}$ and $\de_n V_a = V_n = 0$ on $\de \mathbb{R}^n_+$ for $a = 1, \cdots, n-1,$
    then
    \begin{align}
     \de_n \psi 
      =& \frac{n}{n-2} \psi v_{\gep}^{-1} \de_n v_{\gep} -\frac{1}{n} div V \de_n v_{\gep}
      + \de_n v_{\gep} \de_n V_n -\frac{n}{2(n-1)} \de_n v_{\gep} S_{nn} \notag\\
    =& \frac{n}{n-2} \de_n v_{\gep} v_{\gep}^{-1} \psi -\frac{1}{2(n-1)} \de_n v_{\gep} S_{nn}. 
    \notag
    \end{align}
  \epf
 
  \begin{prop} \label{p:b2}
  Let $\xi_i$ be defined as in Proposition~\ref{p:b1}.
  It follows for $ x \in \de \mathbb{R}^n_+,$
    $$\xi_n(x)= -\frac{n+2}{2 (n-2)} v_{\gep}(x) \de_n v_{\gep}(x)
    S_{nn}(x)^2
    + \frac{4n (n-1)}{(n-2)^2} v_{\gep}(x)^{-1} \de_n v_{\gep}(x) \psi(x)^2.$$
  \end{prop}
  \bpf
    Since $H_{in} = 0$  for $i= 1, \cdots, n$ and $x \in \mathbb R^n_+,$ and    
    $S_{na} = T_{na} = 0$ for  $a= 1, \cdots, n-1$ and $x \in \de \mathbb R^n_+,$ we have
   \beqn
   \xi_n &=& - \frac{1}{2} v_{\gep}^2 \sum_{a, b = 1}^{n-1} \de_n S_{ab} H_{ab} - v_{\gep}
   \psi \de_n S_{nn} + v_{\gep} \de_n \psi S_{nn} + \de_n v_{\gep} \psi
   S_{nn}\\
    &+& \frac{1}{4} v_{\gep}^2 (\sum_{a, b = 1}^{n-1} \de_n S_{ab} S_{ab} + \de_n S_{nn}
    S_{nn}) -\frac{1}{2} v_{\gep}^2 \de_n S_{nn} S_{nn} - v_{\gep}
    \de_n v_{\gep} S_{nn} S_{nn}\\
    &-&  \frac{2}{n-2} v_{\gep} \de_n v_{\gep} S_{nn} S_{nn}
    +  \frac{4 (n-1)}{n-2} (- \de_n v_{\gep} \psi S_{nn} + \psi \de_n \psi).
   \eeqn
   By $\de_n S_{nn}= - \frac{2n}{n-2} v_{\gep}^{-1} \de_n v_{\gep}
  S_{nn}$ and $\de_n S_{a b} = - \frac{2n}{n-1} v_{\gep}^{\frac{2}{n-2}} S_{nn}
   \gd_{ab},$ we get
   \beqn
   \xi_n &=&  \frac{n}{n-1} v_{\gep}^2 v_{\gep}^{\frac{2}{n-2}} S_{nn}\sum_{a= 1}^{n-1} H_{aa}
   + \frac{2n}{n-2}  \psi \de_n v_{\gep} S_{nn} +
    v_{\gep} \de_n \psi S_{nn} + \de_n v_{\gep} \psi
   S_{nn}\\
    &-& \frac{1}{2} v_{\gep}^2 (\frac{n}{n-1} v_{\gep}^{\frac{2}{n-2}} S_{nn} \sum_{a= 1}^{n-1} S_{aa} + \frac{n}{n-2} v_{\gep}^{-1} \de_n v_{\gep}
    S_{nn}^2 ) +\frac{n}{n-2} v_{\gep} \de_n v_{\gep} S_{nn}^2
    \\
    &-&  v_{\gep}
    \de_n v_{\gep} S_{nn} S_{nn}+  \frac{4 (n-1)}{n-2} (- \de_n v_{\gep} \psi S_{nn} + \psi \de_n \psi)
    - \frac{2}{n-2} v_{\gep} \de_n v_{\gep} S_{nn} S_{nn}.
  \eeqn
  Thus,
  \beqn
   \xi_n &=&  \frac{n}{n-1} v_{\gep}^2 v_{\gep}^{\frac{2}{n-2}} S_{nn}\sum_{a= 1}^{n-1} H_{aa}
    - \psi \de_n v_{\gep} S_{nn} +    v_{\gep} \de_n \psi S_{nn}\\
    &-& \frac{1}{2} v_{\gep}^2 (\frac{n}{n-1} v_{\gep}^{\frac{2}{n-2}} S_{nn} \sum_{a= 1}^{n-1} S_{aa} + \frac{n}{n-2} v_{\gep}^{-1} \de_n v_{\gep}
    S_{nn}^2 )+ \frac{4 (n-1)}{n-2} \psi \de_n \psi.
   \eeqn
  By $\sum_{a= 1}^{n-1} H_{aa}= \sum_{i= 1}^n S_{ii} =0$ and
  (\ref{e:v_bdy}),
   we get
   \beqn
   \xi_n
    &=& - \de_n v_{\gep} \psi S_{nn} +  v_{\gep} \de_n \psi S_{nn}-\frac{1}{2} v_{\gep}  (\frac{n}{(n-1)(n-2)}  \de_n v_{\gep} S_{nn}^2
    +\frac{n}{n-2}  \de_n v_{\gep} S_{nn}^2)\\
    &+& \frac{4(n-1)}{n-2} \psi \de_n
    \psi\\
    &=& -\de_n v_{\gep} \psi S_{nn}+ v_{\gep} \de_n
    \psi S_{nn}  - \frac{n^2}{2(n-1)(n-2)}
     v_{\gep} \de_n v_{\gep} S_{nn}^2  + \frac{4(n-1)}{n-2} \psi \de_n \psi.
   \eeqn
    Finally, by
 $\de_n \psi = -\frac{1}{2 (n-1)} \de_n v_{\gep} S_{nn} + \frac{n}{n-2} v_{\gep}^{-1}
  \de_n v_{\gep} \psi,$ we arrive at
    \beqn
   \xi_n
    &=& -\de_n v_{\gep} \psi S_{nn}  + v_{\gep} (-\frac{1}{2 (n-1)} \de_n v_{\gep} S_{nn} + \frac{n}{n-2} v_{\gep}^{-1}
  \de_n v_{\gep} \psi) S_{nn}\\
  & -& \frac{n^2}{2(n-1)(n-2)}
     v_{\gep} \de_n v_{\gep} S_{nn}^2+ \frac{4(n-1)}{n-2} \psi (-\frac{1}{2 (n-1)} \de_n v_{\gep} S_{nn}
      + \frac{n}{n-2} v_{\gep}^{-1}
  \de_n v_{\gep} \psi)\\
  &=& - \frac{n+2}{2(n-2)} v_{\gep} \de_n v_{\gep} S_{nn}^2 +
  \frac{4n (n-1)}{(n-2)^2} v_{\gep}^{-1} \de_n v_{\gep} \psi^2.
   \eeqn
   \epf
 \section{Main estimates} \label{s:est}
  In this section, we assume $g$ is the metric in conformal Fermi
  coordinates as described in Section~\ref{s:construct}. Suppose $V$
  is a smooth vector field which satisfies (\ref{e:star}) and (\ref{i:V}).
  We adopt the notation in Section~\ref{s:construct}.

 \begin{prop}  \label{p:local} There exist positive numbers $\theta,
 C$ and $\gd_0$ such that
    \begin{align}
    &  \int_{B_{\gd} \cap \mathbb{R}^n_+} (\frac{4(n-1)}{n-2} |d (v_{\gep}+ \psi)|_g^2 + R_g (v_{\gep}+ \psi)^2) dx \notag\\
    \leq& \;
      4(n-1)
   \int_{B_{\gd}\cap \de \mathbb{R}^n_+} v_{\gep}^{\frac{2}{n-2}} (v_{\gep}^2 + 2 v_{\gep} \psi + \frac{n}{n-2} \psi^2
   - \frac{n-2}{8 (n-1)^2} v_{\gep}^2 |S_{nn}|^2) d\gs \notag\\
   +& \;\int_{\de B_{\gd} \cap \mathbb R^n_+} \sum_{i=1}^n( \frac{4 (n-1)}{n-2} v_{\gep} \de_i v_{\gep}
    + v_{\gep}^2 \de_k h_{ik}- \de_k v_{\gep}^2 h_{ik}) \frac{x_i}{|x|}
     d \sigma \notag\\
  -& \; \theta \sum_{i, k=1}^n \sum_{2 \leq |\ga| \leq d} |h_{ik, \ga}|^2 \gep^{n-2} \int_{B_{\gd} \cap \mathbb{R}^n_+} (\gep+ |x|)^{2 |\ga| + 2 -2n} dx \notag\\
  +& \; C \gep^{n-2} \sum_{i, k=1}^n \sum_{2 \leq |\ga|\leq d}  |h_{ik, \ga}| \gd^{|\ga|+2-n} + C \gep^{n-2} \gd^{2d+4-n} \notag
 \end{align} for $0 < 2 \epsilon \leq \gd \leq \gd_0,$ where $\theta= \theta
 (n)$, $C=C(n, g)$ and $\gd_0 = \gd_0 (n, g).$
  \end{prop}
  \bpf We write
  $$\frac{4(n-1)}{n-2} |d (v_{\gep}+ \psi)|_g^2 + R_g (v_{\gep}+
  \psi)^2= \frac{4 (n-1)}{n-2} |\de v_{\gep}|^2 + J_1 + J_2 + J_3 + J_4,$$
  where
  \begin{align}
   J_1 =&  \frac{8 (n-1)}{n-2} \sum_{i= 1}^n\de_i v_{\gep} \de_i \psi + \sum_{i,k=1}^n(-\frac{4(n-1)}{n-2} \de_i v_{\gep} \de_k  h_{ik}
   + v_{\gep}^2  \de_i \de_k h_{ik}) \notag\\
   & - \sum_{i,k,l =1}^n (v_{\gep}^2 \de_k (H_{ik} \de_l
   H_{il})
    +\de_k v_{\gep}^2 H_{ik} \de_l H_{il}), \notag
  \end{align}
   \begin{align}
   J_2 = & \sum_{i,k, l = 1}^n (-\frac{1}{4} v_{\gep}^2 \de_l H_{ik} \de_l H_{ik} + \frac{1}{2} v_{\gep}^2 \de_k H_{ik} \de_l H_{il} + \de_k v_{\gep}^2 H_{ik} \de_l H_{il}
           + \frac{2(n-1)}{n-2}  \de_k v_{\epsilon} \de_l v_{\epsilon} H_{ik} H_{il})\notag \\
      &+  \sum_{i,k = 1}^n (  2 v_{\gep} \psi \de_i \de_k H_{ik} -  \frac{8 (n-1)}{n-2} \de_i v_{\gep} \de_k \psi  H_{ik}) + \frac{4(n-1)}{n-2}
      |d\psi|^2, \notag
   \end{align}
   \begin{align}
    J_3 =& \frac{4(n-1)}{n-2} \sum_{i,k= 1}^n(g^{ik} - \gd_{ik}+ h_{ik} -\frac{1}{2}\sum_{l=1}^n H_{il} H_{kl}) \de_i v_{\gep} \de_k
    v_{\gep}\notag\\
       &+(R_g  -\sum_{i,k= 1}^n \de_i \de_k h_{ik} + \sum_{i,k,l= 1}^n\de_k (H_{ik} \de_l H_{il})- \frac{1}{2} (div H)^2 + \frac{1}{4} |\de H|^2)
       v_{\gep}^2,\notag
   \end{align} and
   \begin{align}
    J_4 =&  \frac{8(n-1)}{n-2} \sum_{i,k = 1}^n(g^{ik} - \gd_{ik} + H_{ik}) \de_i v_{\gep} \de_k
    \psi
      + 2 (R_g - \sum_{i,k =1}^n\de_i \de_k H_{ik}) v_{\gep} \psi \notag\\
     &  + R_g \psi^2 + \frac{4(n-1)}{n-2}\sum_{i,k = 1}^n (g^{ik}- \gd_{ik}) \de_i
     \psi
     \de_k \psi. \notag
    \end{align}

   We compute
   \begin{align}
   J_1 =&  \frac{8 (n-1)}{n-2} \sum_{i= 1}^n \de_i (\de_i v_{\gep} \psi)
   -\frac{8(n-1)}{n-2} \Delta v_{\gep} \psi+
   \sum_{i,k=1}^n
    (\de_i (v_{\gep}^2 \de_k h_{ik}) -\de_k (\de_i v_{\gep}^2 h_{ik})) \notag\\
   & +2 \sum_{i,k = 1}^n(v_{\gep} \de_i \de_k v_{\gep} -\frac{n}{n-2} \de_i v_{\gep} \de_k
   v_{\gep}) h_{ik}
   - \sum_{i,k,l =1}^n \de_k (v_{\gep}^2  H_{ik} \de_l
   H_{il}). \notag
   \end{align}
  By (\ref{e:v}) and (\ref{e:hess_v}),
  \begin{align}
   J_1 \leq &   \frac{8 (n-1)}{n-2} \sum_{i= 1}^n \de_i (\de_i v_{\gep}
   \psi) +   \sum_{i,k=1}^n
    (\de_i (v_{\gep}^2 \de_k h_{ik}) -\de_k (\de_i v_{\gep}^2
    h_{ik})) - \sum_{i,k,l =1}^n \de_k (v_{\gep}^2 H_{ik} \de_l
   H_{il})\notag \\
   &  + C \gep^{n-2} (\gep + |x|)^{2d +4 -2n}. \notag
  \end{align}
  Thus, integrating $J_1$ over $B_{\gd} \cap \mathbb{R}^n_+$ and using (\ref{i:V}),
  \beqn
  \int_{B_{\gd} \cap \mathbb{R}^n_+} J_1 dx & \leq &
   \int_{\de B_{\gd} \cap \mathbb{R}^n_+} \sum_{i,k = 1}^n(v_{\gep}^2 \de_k h_{ik} -\de_k v_{\gep}^2 h_{ik})
   \frac{x_i}{|x|} d\gs- \int_{B_{\gd} \cap \de \mathbb{R}^n_+}
   \frac{8 (n-1)}{n-2} \de_n v_{\gep} \psi d\gs\\
   &+& C \sum_{i,k = 1}^n \sum_{2 \leq |\ga| \leq d} |h_{ik, \ga}|\gd^{|\ga|+2 -n} \gep^{n-2} + C \gd^{2d + 4
   -n} \gep^{n-2}.
  \eeqn

  For $J_2,$ we first note that by Proposition~\ref{p:b1} and (\ref{e:T}), $J_2 = - \frac{1}{4} |P|^2 + \sum_{i= 1}^n \de_i \xi_i.$
  And by (\ref{i:V})
 $$\int_{B_{\gd} \cap \mathbb{R}^n_+} \xi_i \frac{x_i}{|x|} d \gs
 \leq C \sum_{i,k =1}^n \sum_{2 \leq |\ga| \leq d} |h_{ik, \ga}|^2 \gd^{2 |\ga|+ 2- n}
 \gep^{n-2}.$$ Moreover,  by Proposition~\ref{p:P} there exists $\theta> 0$  such that
 $$8 \theta \sum_{i, j}  \sum_{2 \leq |\ga| \leq d} |h_{ij, \ga}|^2 \gep^{n-2} \int_{B_{\gd} \cap \mathbb{R}^n_+} (\gep + |x|) ^{2 |\ga| +2 -2n} dx  \leq
      \int_{B_{\gd} \cap \mathbb{R}^n_+}  |P|^2 dx.$$
 Hence, using Proposition~\ref{p:b2}
  \beqn
  \int_{B_{\gd} \cap \mathbb{R}^n_+} J_2 dx &=& -\int_{B_{\gd} \cap \mathbb{R}^n_+}\frac{1}{4} |P|^2 dx +
  \int_{\de B_{\gd} \cap \mathbb{R}^n_+} \xi_i \frac{x_i}{|x|} d\gs
  - \int_{B_{\gd} \cap \de \mathbb{R}^n_+} \xi_n d\gs\\
  &\leq&
  \int_{B_{\gd}\cap \de \mathbb{R}^n_+} (\frac{n+2}{2 (n-2)}v_{\gep} \de_n v_{\gep} |S_{nn}|^2
  - \frac{4n (n-1)}{(n-2)^2} v_{\gep}^{-1} \de_n v_{\gep} \psi^2)
  d\gs\\
  & &-\, 2 \theta \sum_{i, k=1}^n \sum_{2 \leq |\ga| \leq d} |h_{ik, \ga}|^2 \gep^{n-2} \int_{B_{\gd} \cap \mathbb{R}^n_+} (\gep+ |x|)^{2 |\ga| + 2 -2n}
  dx\\
  & &+ \,C \sum_{i,k = 1}^n \sum_{2 \leq |\ga| \leq d} |h_{ik, \ga}|^2 \gd^{2 |\ga| + 2
  -n} \gep^{n-2}.
  \eeqn

  For $J_3$ and $J_4,$ by (\ref{i:V}), Proposition~\ref{p:scalar} and Cauchy inequality,
   \begin{align}
   J_3 + J_4 &\leq C \sum_{i,k = 1}^n \sum_{2 \leq |\ga| \leq d} (|h_{ik, \ga}|^2 (\gep + |x|)^{2 |\ga|
   + 4 -2 n}  +  |h_{ik, \ga}| (\gep + |x|)^{|\ga|+ d+
   3 -2n}) \gep^{n-2} \notag\\
   &+ C (\gep + |x|)^{2d + 4 -2n} \gep^{n-2}\notag \\
    & \leq \theta \sum_{i, k=1}^n \sum_{2 \leq |\ga| \leq d} |h_{ik, \ga}|^2 \gep^{n-2} (\gep+ |x|)^{2 |\ga| + 2
    -2n}
    + C (\gep + |x|)^{2d + 4 -2n} \gep^{n-2}.\notag
   \end{align}
   Thus,
   $$\int_{B_{\gd} \cap \mathbb{R}^n_+} (J_3 + J_4) dx  \leq \theta \sum_{i, k=1}^n \sum_{2 \leq |\ga| \leq d} |h_{ik, \ga}|^2 \gep^{n-2} \int_{B_{\gd} \cap \mathbb{R}^n_+} (\gep+ |x|)^{2 |\ga| + 2
    -2n} dx + C \gd^{2d + 4 -n} \gep^{n-2}.$$

   Finally, by (\ref{e:v}) we compute
   $$\int_{B_{\gd} \cap \mathbb{R}^n_+} \frac{4 (n-1)}{n-2} |d v_{\gep}|^2 dx = \frac{4(n-1)}{n-2}
   (\int_{B_{\gd}\cap \de \mathbb{R}^n_+} -v_{\gep} \de_n v_{\gep} d\gs +
    \int_{\de B_{\gd} \cap \mathbb R^n_+} \sum_{i= 1}^n v_{\gep}  \de_i v_{\gep} \frac{x_i}{|x|} d\gs).$$

   Combining the above, we obtain
   \begin{align}
    &  \int_{B_{\gd} \cap \mathbb{R}^n_+} (\frac{4(n-1)}{n-2} |d (v_{\gep}+ \psi)|_g^2 + R_g (v_{\gep}+ \psi)^2) dx \notag\\
    \leq & \;- \frac{4(n-1)}{n-2} \int_{B_{\gd}\cap \de \mathbb{R}^n_+}
   ( v_{\gep} \de_n v_{\gep} + 2 \de_n v_{\gep} \psi
   + \frac{n}{(n-2)} v_{\gep}^{-1} \de_n v_{\gep} \psi^2) d\gs\notag\\
    + & \;  \frac{n+2}{2 (n-1)}\int_{B_{\gd}\cap \de \mathbb{R}^n_+} v_{\gep} \de_n v_{\gep}  |S_{nn}|^2
    d\gs\notag\\
   +& \int_{\de B_{\gd} \cap \mathbb R^n_+} \sum_{i=1}^n  ( \frac{4 (n-1)}{n-2} v_{\gep} \de_i
   v_{\gep}+ v_{\gep}^2 \de_k h_{ik}- \de_k v_{\gep}^2 h_{ik}) \frac{x_i}{|x|}
     d \sigma
   \notag\\
  -& \;
     \theta \sum_{i, k=1}^n \sum_{2 \leq |\ga| \leq d} |h_{ik, \ga}|^2 \gep^{n-2} \int_{B_{\gd} \cap \mathbb{R}^n_+} (\gep+ |x|)^{2 |\ga| + 2 -2n}
     dx\notag\\
  +& \;C \gep^{n-2} \sum_{i, k=1}^n \sum_{2 \leq |\ga|\leq d}  |h_{ik, \ga}| \gd^{-n +2+ |\ga|}
  + C \gep^{n-2} \gd^{2d+4-n}. \notag
 \end{align}
  Finally, by (\ref{e:v_bdy}) and $v_{\gep} \de_n v_{\gep} |S_{nn}|^2 \leq 0$
  for  $x \in \de \mathbb{R}^n_+,$
  \begin{align}
   &  - \frac{4(n-1)}{n-2} \int_{B_{\gd}\cap \de \mathbb{R}^n_+}
   ( v_{\gep} \de_n v_{\gep} + 2 \de_n v_{\gep} \psi
   + \frac{n}{(n-2)} v_{\gep}^{-1} \de_n v_{\gep} \psi^2) d\gs\notag\\
   & +    \frac{n+2}{2 (n-1)}\int_{B_{\gd}\cap \de \mathbb{R}^n_+} v_{\gep} \de_n v_{\gep}  |S_{nn}|^2
    d\gs \notag\\
  & \leq - \frac{4(n-1)}{n-2} \int_{B_{\gd}\cap \de \mathbb{R}^n_+}
   ( v_{\gep} \de_n v_{\gep} + 2 \de_n v_{\gep} \psi
   + \frac{n}{(n-2)} v_{\gep}^{-1} \de_n v_{\gep} \psi^2) d\gs\notag\\
   & +    \frac{1}{2 (n-1)}\int_{B_{\gd}\cap \de \mathbb{R}^n_+} v_{\gep} \de_n v_{\gep}  |S_{nn}|^2
    d\gs \notag\\
   &= 4(n-1)
   \int_{B_{\gd}\cap \de \mathbb{R}^n_+} v_{\gep}^{\frac{2}{n-2}} (v_{\gep}^2 + 2 v_{\gep} \psi + \frac{n}{n-2} \psi^2
   - \frac{n-2}{8 (n-1)^2} v_{\gep}^2 |S_{nn}|^2) d\gs.\notag
  \end{align} This completes the proof.
  \epf

\begin{prop}\label{p:local2}
 \begin{align}
 & 4(n-1) \int_{B_{\gd}\cap \de \mathbb{R}^n_+} v_{\gep}^{\frac{2}{n-2}} (v_{\gep}^2 + 2 v_{\gep} \psi + \frac{n}{n-2} \psi^2
   - \frac{n-2}{8 (n-1)^2} v_{\gep}^2 S_{nn}^2) d\gs \notag\\
 & \leq \mathcal{\mathcal{Q}} (B, \de B) (\int_{B_{\gd} \cap \de \mathbb{R}^n_+}
   (v_{\gep} + \psi)^{\frac{2(n-1)}{n-2}} d\gs)^{\frac{n-2}{n-1}} + C
   \sum_{i,k= 1}^n \sum_{2 \leq |\ga| \leq d} |h_{ik, \ga}| \gd^{|\ga| - n
   +1} \gep^{n-1} \notag\\
 & + C\sum_{i, k=1}^n \sum_{2 \leq |\ga| \leq d} |h_{ik, \ga}|^2
 \gep^{n-1} \gd^2
      \int_{B_{\gd} \cap \de \mathbb{R}^n_+} (\gep+ |x|)^{2 |\ga|-2n+ 2}
      d\gs \notag
 \end{align} for $0 < 2 \gep \leq \gd \leq \gd_0$ and $\gd_0$
 sufficiently small.
\end{prop}
 \bpf
 Recall that  $$\mathcal{Q}(B^n, \de B^n) \left( \int_{\de \mathbb R^n_+}  v_{\gep}^{\frac{2(n-1)}{n-2}} d \gs\right)^{\frac{n-2}{n-1}} = \frac{4(n-1)}{n-2}\int_{\mathbb R^n_+} |\gra v_{\gep}|^2 dx$$ and
   $$\int_{\mathbb R^n_+} |\gra v_{\gep}|^2 dx = (n-2) \int_{\de \mathbb R^n_+} v_{\gep}^{\frac{2(n-1)}{n-2}} d\gs.$$
  Then it follows that
  $$4 (n-1) ( \int_{\de \mathbb{R}^n_+} v_{\gep}^{\frac{2(n-1)}{n-2}} d\gs)^{\frac{1}{n-1}} = \mathcal{Q}(B, \de
  B).$$  Besides,  
  since $V_n = 0$ on $\de \mathbb{R}^n_+,$
  we have
  \begin{align}
  \psi &=  \frac{n-2}{2(n-1)} v_{\gep}^{-\frac{n}{n-2}}
  \sum_{a= 1}^{n-1}\de_a (v_{\gep}^{\frac{2(n-1)}{n-2}} V_a)+ \frac{n-2}{4(n-1)} v_{\gep} (2 \de_n V_n - \frac{2}{n} div V)\notag\\ 
  &=  \frac{n-2}{2(n-1)} v_{\gep}^{-\frac{n}{n-2}}
  \sum_{a= 1}^{n-1}\de_a (v_{\gep}^{\frac{2(n-1)}{n-2}} V_a)+ \frac{n-2}{4(n-1)} v_{\gep} S_{nn} \notag
  \end{align}
   for $x \in \de \mathbb{R}^n_+.$
   Moreover, by (\ref{i:V})
   $$\int_{\de B_{\gd} \cap \de \mathbb{R}^n_+} v_{\gep}^{\frac{2(n-1)}{n-2}}\sum_{a=1}^{n-1} V_a \frac{x_a}{|x|} d \mu
   \leq C \sum_{i,k= 1}^n \sum_{2 \leq |\ga| \leq d} |h_{ik, \ga}| \gd^{|\ga| - n
   +1} \gep^{n-1}.$$
   Thus,
   $$\int_{ B_{\gd} \cap \de \mathbb{R}^n_+} 2  v_{\gep}^{\frac{n}{n-2}}  \psi d \gs
   - \int_{ B_{\gd} \cap \de \mathbb{R}^n_+} \frac{n-2}{2(n-1)} v_{\gep}^{\frac{2(n-1)}{n-2}}  S_{nn} d \gs
   \leq  C \sum_{i,k= 1}^n \sum_{2 \leq |\ga| \leq d} |h_{ik, \ga}| \gd^{|\ga| - n
   +1} \gep^{n-1}.$$
   Putting above together and using Holder inequality, we get
   \beqn
   & & 4(n-1) \int_{B_{\gd}\cap \de \mathbb{R}^n_+} v_{\gep}^{\frac{2}{n-2}} (v_{\gep}^2 + 2 v_{\gep} \psi + \frac{n}{n-2} \psi^2
   - \frac{n-2}{8 (n-1)^2} v_{\gep}^2 S_{nn}^2) d\gs\\
    &\leq& 4(n-1) \int_{B_{\gd}\cap \de \mathbb{R}^n_+} v_{\gep}^{\frac{2}{n-2}} (v_{\gep}^2 +  \frac{n-2}{2(n-1)} v_{\gep}^2 S_{nn}+ \frac{n}{n-2} \psi^2
   - \frac{n-2}{8 (n-1)^2} v_{\gep}^2 S_{nn}^2) d\gs\\
   & & +C \sum_{i,k= 1}^n \sum_{2 \leq |\ga| \leq d} |h_{ik, \ga}| \gd^{|\ga| - n
   +1} \gep^{n-1}\\
   &\leq& \mathcal{Q}(B, \de B) \left( \int_{B_{\gd}\cap \de \mathbb{R}^n_+}(v_{\gep}^2
   + \frac{n-2}{2(n-1)} v_{\gep}^2 S_{nn} + \frac{n}{n-2} \psi^2
   - \frac{n-2}{8 (n-1)^2} v_{\gep}^2 S_{nn}^2)^{\frac{n-1}{n-2}}
   d\gs \right)^{\frac{n-2}{n-1}}\\
   & & +C \sum_{i,k= 1}^n \sum_{2 \leq |\ga| \leq d} |h_{ik, \ga}| \gd^{|\ga| - n
   +1} \gep^{n-1}.
   \eeqn

    We next notice that by Taylor expansion, there exists a constant $C_0= C_0(n)$ such that
    \begin{align}
    &(1 +\frac{n-2}{2(n-1)} y + \frac{n}{n-2} z^2 - \frac{n-2}{8(n-1)^2} y^2)^{\frac{n-1}{n-2}}
     -(1+ z)^{\frac{2(n-1)}{n-2}} + \frac{2(n-1)}{n-2} z -\frac{1}{2} y \notag\\
      \leq & \;C_0 (|y|^3 + |z|^3) \notag
    \end{align} for $|y|, |z| \leq \frac{1}{2}.$
   By (\ref{i:V}), $|S_{nn}| \leq \frac{1}{2}$ and $|\psi| \leq \frac{1}{2} v_{\gep}$
   for $|x| \leq \gd.$ Hence,
  \begin{align}
   & (v_{\gep}^2
   + \frac{n-2}{2(n-1)} v_{\gep}^2 S_{nn} + \frac{n}{n-2} \psi^2
   - \frac{n-2}{8 (n-1)^2} v_{\gep}^2 S_{nn}^2)^{\frac{n-1}{n-2}} -
   (v_{\gep} + \psi)^{\frac{2(n-1)}{n-2}}\notag\\
    & \;+  \frac{2(n-1)}{n-2}  v_{\gep}^{\frac{n}{n-2}}  \psi
    - \frac{1}{2} v_{\gep}^{\frac{2(n-1)}{n-2}} S_{nn} \notag\\
   \leq &\; C_0 v_{\gep}^{\frac{2(n-1)}{n-2}} (|\frac{\psi}{v_{\gep}}|^3 +
   |S_{nn}|^3)
   \leq  C v_{\gep}^{\frac{2(n-1)}{n-2}} (|\frac{\psi}{v_{\gep}}|^2 +
   |S_{nn}|^2) \gd^2 \notag\\
   \leq & \;
    C \sum_{i,k= 1}^n \sum_{2 \leq |\ga| \leq d} |h_{ik, \ga}|^2 (\gep + |x|)^{2|\ga| - 2n
   +2} \gep^{n-1} \gd^2. \notag
   \end{align}
 Thus,
   \beqn
   & & \int_{B_{\gd}\cap \de \mathbb{R}^n_+}(v_{\gep}^2
   + \frac{n-2}{2(n-1)} v_{\gep}^2 S_{nn} + \frac{n}{n-2} \psi^2
   - \frac{n-2}{8 (n-1)^2} v_{\gep}^2 S_{nn}^2)^{\frac{n-1}{n-2}}
   d\gs\\
   & \leq&  \int_{B_{\gd}\cap \de \mathbb{R}^n_+} (v_{\gep} +
   \psi)^{\frac{2(n-1)}{n-2}} d\gs   
    +  \int_{B_{\gd}\cap \de \mathbb{R}^n_+}  \frac{2(n-1)}{n-2}  v_{\gep}^{\frac{n}{n-2}}  \psi d \gs
    -  \int_{B_{\gd}\cap \de \mathbb{R}^n_+} \frac{1}{2} v_{\gep}^{\frac{2(n-1)}{n-2}} S_{nn} d \gs\\ 
   &+& C \sum_{i,k= 1}^n \sum_{2 \leq |\ga| \leq d} |h_{ik, \ga}|^2  \gep^{n-1} \gd^2 \int_{B_{\gd} \cap \de \mathbb{R}^n_+} (\gep + |x|)^{2|\ga| - 2n
   +2} d \gs\\   
   & \leq& \int_{B_{\gd}\cap \de \mathbb{R}^n_+} (v_{\gep} +
   \psi)^{\frac{2(n-1)}{n-2}} d\gs +C \sum_{i,k= 1}^n \sum_{2 \leq |\ga| \leq d} |h_{ik, \ga}|^2  \gep^{n-1} \gd^2 \int_{B_{\gd} \cap \de \mathbb{R}^n_+} (\gep + |x|)^{2|\ga| - 2n
   +2} d \gs\\
   &+&   C \sum_{i,k= 1}^n \sum_{2 \leq |\ga| \leq d} |h_{ik, \ga}| \gd^{|\ga| - n
   +1} \gep^{n-1}.    \eeqn
   This completes the proof.
  \epf


\section{Proof of the main theorems} \label{s:proof}
  In this section, we construct a test function $\phi_{(\gep, \gd)}$
  with energy functional less than $\mathcal{Q}(B, \de B)$ and prove 
  Theorem~\ref{t:main} and \ref{t:main2}.
  Since the case that $\mathcal{Q} (M, \de M, g) \leq 0$ 
 is trivial,   it suffices to
 consider $\mathcal{Q} (M, \de M, g) > 0.$
 
  After a conformal change of the metric, we may assume $\de M$ is
  totally geodesic. Let $p \in \de M$ and let $(x_1, \cdots, x_n)$
  be the conformal Fermi coordinates around $p$ described in
  Section~\ref{s:construct}.
  We denote by $G$ the Green's function
   of the  conformal Laplacian with pole at $p$ which satisfies the Neumann boundary condition. We assume
  that $G$ is normalized such that $\lim_{|x| \rightarrow 0} |x|^{n-2} G (x) = 1.$
  Then $G$ satisfies \cite{Bchen09}
 \beq \label{i:Gexp}
  |G(x)- |x|^{2-n}| \leq C \sum_{i,k= 1}^n \sum_{2 \leq |\ga| \leq
  d} |h_{ik, \ga}| |x|^{|\ga|+2 -n} + C |x|^{d+3-n}.
 \eeq
 Moreover, we define as in \cite{Bchen09} a flux integral
 \beqn
  \mathcal{I}(p, \gd) &=& \frac{4(n-1)}{n-2}\int_{\de B_{\gd} \cap \mathbb{R}^n_+}
     \sum_{i=1}^n(|x|^{2-n} \de_i G - G \de_i |x|^{2-n}) \frac{x_i}{|x|}
     d\gs\\
     & & - \int_{\de B_{\gd} \cap \mathbb{R}^n_+}
     \sum_{i, k=1}^n|x|^{2-2n} (|x|^2 \de_k h_{ik} - 2n x_k h_{ik})
     \frac{x_i}{|x|} d\gs
 \eeqn for $\gd > 0$ sufficiently small.

  We define
  $$\phi_{(\gep, \gd)} = \eta_{\gd} (v_{\gep} + \psi) + (1- \eta_{\gd}) \gep^{\frac{n-2}{2}}
  G,$$ where $\psi$ is the function constructed in Section~\ref{s:construct}.
   We recall that $$\gep^{\frac{n-2}{2}} (\gep + |x|)^{-n+2} \leq v_{\gep} (x) \leq C(n) \gep^{\frac{n-2}{2}} (\gep + |x|)^{-n+2}
   \quad \text{for}\; x \in \mathbb{R}^n_+;$$
   $$|\de v_{\gep}| (x) \leq C(n) \gep^{\frac{n-2}{2}} (\gep +
   |x|)^{-n+1} \quad \text{for}\; x \in \mathbb{R}^n_+;$$
   and 
   \beq \label{i:v_exp}
   |v_{\gep} - \gep^{\frac{n-2}{2}} |x|^{-n+2}| \leq C(n) \gep^{\frac{n}{2}} |x|^{-n+1}
    \quad \text{for}\;  x \in \mathbb{R}^n_+,\; \text{and}\; |x|  \geq 2\gep.\eeq
\begin{prop}  \label{p:main}
    \begin{align}
    &  \int_M (\frac{4(n-1)}{n-2} |d \phi_{(\gep, \gd)}|_g^2 + R_g \phi_{(\gep, \gd)}^2) dV_g \notag\\
    \leq& \;
   \mathcal{\mathcal{Q}} (B, \de B) (\int_{\de M}
   \phi_{(\gep, \gd)}^{\frac{2(n-1)}{n-2}}d\gs_g)^{\frac{n-2}{n-1}}
      -  \frac{\theta}{2} \sum_{i, k=1}^n \sum_{2 \leq |\ga| \leq d} |h_{ik, \ga}|^2 \gep^{n-2}
      \int_{B_{\gd} \cap \mathbb{R}^n_+} (\gep+ |x|)^{2 |\ga| + 2 -2n} dx
     \notag\\
  - & \; \gep^{n-2} \mathcal{I} (p, \gd)+
      C \gep^{n-2} \sum_{i, k=1}^n \sum_{2 \leq |\ga|\leq d}  |h_{ik, \ga}| \gd^{-n +2+ |\ga|} + C \gep^{n-2} \gd^{2d+4-n}
  + C \gd^{-n+1} \gep^{n-1} \notag
 \end{align} for $0 < 2\epsilon \leq \gd\leq \gd_0$ and $\gd_0$
 sufficiently small.
  \end{prop}
 \bpf
  Let $\Omega_{\gd}$ be the coordinates ball of radius $\gd$ in Fermi
  coordinates.
  In other words, $(x_1, \cdots, x_n)$
  satisfies $x_1^2 + \cdots + x_n^2 < \gd^2$ and $x_n \geq 0.$
  By divergence theorem
   \beqn
  & & \int_{M \setminus \Omega_{\gd}} (\frac{4(n-1)}{n-2} |\gra_g \phi_{(\gep, \gd)}|^2 + R_g \phi_{(\gep, \gd)}^2)
  dV_g\\
    & =&   \int_{M \setminus \Omega_{\gd}} -(\frac{4(n-1)}{n-2} \Delta_g \phi_{(\gep, \gd)} - R_g \phi_{(\gep, \gd)}) (\phi_{(\gep, \gd)} - \gep^{\frac{n-2}{2}} G)dV_g
     \\
    & &
    + \frac{4(n-1)}{n-2}\int_{\de (M \setminus \Omega_{\gd})}  (\gra_{\nu_g} \phi_{(\gep, \gd)} \phi_{(\gep, \gd)}
    + \gep^{\frac{n-2}{2}} (\phi_{(\gep, \gd)} \gra_{\nu_g} G - G \gra_{\nu_g} \phi_{(\gep, \gd)})) d\gs_g,
   \eeqn where $\nu_g$ is the unit outer normal on $ \de (M \setminus \Omega_{\gd})$ with respect to $g.$
  Notice that $$\de (M \setminus \Omega_{\gd}) = (\de M \setminus \Omega_{\gd})\cup (\de \Omega_{\gd} \setminus  \de
  M).$$ We will compute the above integral in several steps.

 We first notice that for $x \in M \setminus \Omega_{\gd},$ we have
  $\phi_{(\gep, \gd)} - \gep^{\frac{n-2}{2}} G =
  \eta_{\gd} (v_{\gep} + \psi - \gep^{\frac{n-2}{2}} G).$
 In particular, $\phi_{(\gep, \gd)} - \gep^{\frac{n-2}{2}} G= 0$ in
 $M \setminus \Omega_{2\gd}.$
 By (\ref{i:Gexp}) and (\ref{i:V}),
  \begin{align}
  & \sup_{M \setminus \Omega_{\gd}}(|\phi_{(\gep, \gd)} - \gep^{\frac{n-2}{2}} G| + \gd^2 |\frac{4(n-1)}{n-2} \Delta_g \phi_{(\gep, \gd)} - R_g \phi_{(\gep,
  \gd)}|)\notag\\
   \leq &\; C \sum_{i, k = 1}^n \sum_{2 \leq |\ga| \leq d} |h_{ik,
   \ga}| \gd^{|\ga| +2 -n} \gep^{\frac{n-2}{2}} + C \gd^{d + 3 -n}
   \gep^{\frac{n-2}{2}}+ C \gd^{-n+1} \gep^{\frac{n}{2}}.\notag
  \end{align}
   Thus,
  \beqn
   & &-\int_{M \setminus \Omega_{\gd}} (\frac{4(n-1)}{n-2} \Delta_g \phi_{(\gep, \gd)} - R_g \phi_{(\gep, \gd)}) (\phi_{(\gep, \gd)} - \gep^{\frac{n-2}{2}}
   G)dV_g\\
   &\leq & C  \sum_{i, k = 1}^n \sum_{2 \leq |\ga| \leq d} |h_{ik,
   \ga}|^2 \gd^{2|\ga| +2 -n} \gep^{n-2} + C \gd^{2d +4 -n}
   \gep^{n-2} + C \gd^{-n} \gep^{n}.
  \eeqn

  We now compute the boundary terms on $\de M \setminus \Omega_{\gd}.$
  Since $\gra_{\nu_g} G = 0$ on $\de M,$ by (\ref{e:v_bdy}),
  Proposition~\ref{p:bdy} and (\ref{i:V})
  $$\sup_{\de M \cap (\Omega_{2 \gd} \setminus \Omega_{\gd})} |\gra_{\nu_g} \phi_{(\gep, \gd)}|
   \leq \sup_{\de M \cap (\Omega_{2 \gd} \setminus \Omega_{\gd})} |\de_n v_{\gep}+ \de_n \psi|
   \leq C \gep^{\frac{n}{2}} \gd^{-n} + C \sum_{i,k = 1}^n \sum_{2 \leq |\ga| \leq d}
    |h_{ik, \ga}| \gd^{|\ga|-n} \gep^{\frac{n}{2}}.$$
  Hence,
  \begin{align}
  & \int_{\de M \setminus \Omega_{\gd}}  (\gra_{\nu_g} \phi_{(\gep, \gd)} \phi_{(\gep, \gd)} + \gep^{\frac{n-2}{2}} (\phi_{(\gep, \gd)} \gra_{\nu_g} G - G \gra_{\nu_g} \phi_{(\gep, \gd)}))
  d\gs_g \notag\\
  & =  \int_{\de M \setminus \Omega_{\gd}}  \gra_{\nu_g} \phi_{(\gep,
\gd)} (\phi_{(\gep, \gd)} - \gep^{\frac{n-2}{2}}  G) d \gs_g \notag\\
  &\leq C \sum_{i, k = 1}^n \sum_{2 \leq |\ga| \leq d} |h_{ik,
   \ga}| \gd^{|\ga| +1 -n} \gep^{n-1} + C \gd^{2d +4 -n}
   \gep^{n-2} + C \gd^{-n} \gep^{n}. \notag
  \end{align}

 We next compute the boundary terms on $\de \Omega_{\gd} \setminus  \de
  M.$
  \beqn
  \int_{\de \Omega_{\gd} \setminus  \de M} \gra_{\nu_g} \phi_{(\gep, \gd)} \phi_{(\gep, \gd)}
  d\gs_g
 &\leq&  \int_{\de B_{\gd} \cap\mathbb{R}^n_+} \sum_{i=1}^n(- \de_i v_{\gep} v_{\gep} + \sum_{k=1}^n v_{\gep}
 \de_k v_{\gep} h_{ik}) \frac{x_i}{|x|} d\gs\\
 &+& C \sum_{i,k=1}^n \sum_{2 \leq |\ga| \leq d} |h_{ik, \ga}|
 \gd^{|\ga| +2 -n} \gep^{n-2} + C \gd^{2d+4-n} \gep^{n-2}.
  \eeqn
  Also,
  \begin{align}
 \int_{\de \Omega_{\gd} \setminus  \de
  M} (\phi_{(\gep, \gd)} \gra_{\nu_g} G - G
\gra_{\nu_g} \phi_{(\gep, \gd)}) d\gs_g
 &\leq -\int_{\de B_{\gd} \cap \mathbb{R}^n_+}\sum_{i=1}^n (v_{\gep} \de_i G - G \de_i v_{\gep}) \frac{x_i}{|x|} d
 \gs \notag\\
  &+ C \sum_{i,k=1}^n \sum_{2 \leq |\ga| \leq d} |h_{ik, \ga}|
 \gd^{|\ga| +2 -n} \gep^{\frac{n-2}{2}} + C \gd^{2d+4-n}
 \gep^{\frac{n-2}{2}}. \notag
  \end{align}
  Combining the above, we obtain
  \beqn
  & & \int_{M \setminus \Omega_{\gd}} (\frac{4(n-1)}{n-2} |d \phi_{(\gep, \gd)}|_g^2 + R_g \phi_{(\gep, \gd)}^2) dV_g
    \\
  &\leq& -\frac{4(n-1)}{n-2}\int_{\de B_{\gd} \cap \mathbb{R}^n_+} \sum_{i=1}^n( \de_i v_{\gep} v_{\gep} - \sum_{k=1}^n v_{\gep}
 \de_k v_{\gep} h_{ik} + \gep^{\frac{n-2}{2}} (v_{\gep} \de_i
G - G \de_i v_{\gep})) \frac{x_i}{|x|} d\gs\\
  & +&  C \sum_{i,k=1}^n \sum_{2 \leq |\ga| \leq d} |h_{ik, \ga}|
 \gd^{|\ga| +2 -n} \gep^{n-2} + C \gd^{2d+4-n} \gep^{n-2}+ C
 \gd^{-n} \gep^{n}.
  \eeqn

 On the other hand, by Proposition~\ref{p:local} and \ref{p:local2},
\begin{align}
    &  \int_{\Omega_{\gd}} (\frac{4(n-1)}{n-2} |d\phi_{(\gep, \gd)}|_g^2 + R_g (\phi_{(\gep, \gd)})^2) dV_g \notag\\
\leq & \mathcal{\mathcal{Q}} (B, \de B) (\int_{\de M}
   \phi_{(\gep, \gd)}^{\frac{2(n-1)}{n-2}} d\gs_g)^{\frac{n-2}{n-1}}
 + \int_{\de B_{\gd} \cap \mathbb R^n_+} \sum_{i=1}^n( \frac{4 (n-1)}{n-2} v_{\gep} \de_i v_{\gep}
    + v_{\gep}^2 \de_k h_{ik}- \de_k v_{\gep}^2 h_{ik}) \frac{x_i}{|x|}
     d \sigma \notag\\
  &- \; \theta \sum_{i, k=1}^n \sum_{2 \leq |\ga| \leq d} |h_{ik, \ga}|^2 \gep^{n-2} \int_{B_{\gd} \cap \mathbb{R}^n_+} (\gep+ |x|)^{2 |\ga| + 2 -2n}
  dx\notag\\
  & + C\sum_{i, k=1}^n \sum_{2 \leq |\ga| \leq d} |h_{ik, \ga}|^2
  \gep^{n-1} \gd^2
      \int_{B_{\gd} \cap \de \mathbb{R}^n_+} (\gep+ |x|)^{2 |\ga|-2n+ 2}
      d\gs\notag\\
  & + \; C  \sum_{i, k=1}^n \sum_{2 \leq |\ga|\leq d}  |h_{ik, \ga}| \gd^{-n +2+ |\ga|} \gep^{n-2} + C  \gd^{2d+4-n} \gep^{n-2}. \notag
 \end{align}

 Adding the above two inequalities, we get
  \begin{align}
    &  \int_M (\frac{4(n-1)}{n-2} |d\phi_{(\gep, \gd)}|_g^2 + R_g (\phi_{(\gep, \gd)})^2) dV_g \notag\\
\leq & \mathcal{\mathcal{Q}} (B, \de B) (\int_{\de M}
   \phi_{(\gep, \gd)}^{\frac{2(n-1)}{n-2}} d\gs_g)^{\frac{n-2}{n-1}}
 + \int_{\de B_{\gd} \cap \mathbb R^n_+} \sum_{i=1}^n(
     v_{\gep}^2 \de_k h_{ik}+ \frac{n}{n-2} \de_k v_{\gep}^2 h_{ik}) \frac{x_i}{|x|}
     d \sigma \notag\\
  & -\frac{4(n-1)}{n-2}\int_{\de B_{\gd} \cap \mathbb{R}^n_+} \sum_{i=1}^n
  \gep^{\frac{n-2}{2}} (v_{\gep} \de_i
G - G \de_i v_{\gep}) \frac{x_i}{|x|} d\gs\notag\\
  & - \; \theta \sum_{i, k=1}^n \sum_{2 \leq |\ga| \leq d} |h_{ik, \ga}|^2 \gep^{n-2} \int_{B_{\gd} \cap \mathbb{R}^n_+} (\gep+ |x|)^{2 |\ga| + 2 -2n} dx \notag\\
  & + C\sum_{i, k=1}^n \sum_{2 \leq |\ga| \leq d} |h_{ik, \ga}|^2
  \gep^{n-1} \gd^2
      \int_{B_{\gd} \cap \de \mathbb{R}^n_+} (\gep+ |x|)^{2 |\ga|-2n+ 2}
      d\gs\notag \\
  & + \; C  \sum_{i, k=1}^n \sum_{2 \leq |\ga|\leq d}  |h_{ik, \ga}| \gd^{-n +2+ |\ga|} \gep^{n-2} + C  \gd^{2d+4-n} \gep^{n-2} + C \gd^{-n} \gep^n. \notag
 \end{align}
 
 Since $$\gep^{n-1} \gd^2 \int_{B_{\gd} \cap \de \mathbb{R}^n_+} (\gep+ |x|)^{2 |\ga|-2n+ 2}
      d\gs = \gep^{2 |\ga|} \gd^2 \int_0^{\frac{\gd}{\gep}} (1+t)^{2 |\ga|-2n +2} t^{n-2} dt$$
    and
    $$\gep^{n-2} \int_{B_{\gd} \cap \mathbb{R}^n_+} (\gep+ |x|)^{2 |\ga| + 2 -2n} dx
    = \gep^{2 |\ga|} \int_0^{\frac{\gd}{\gep}} (1+t)^{2 |\ga| -2n +2} t^{n-1} dt,$$
 then for $\gd$ sufficiently small and $2 \gep \leq \gd,$ we have
  $$C \gep^{n-1} \gd^2 \int_{B_{\gd} \cap \de \mathbb{R}^n_+} (\gep+ |x|)^{2 |\ga|-2n+ 2}
      d\gs < \frac{\theta}{2} \, \gep^{n-2} \int_{B_{\gd} \cap \mathbb{R}^n_+} (\gep+ |x|)^{2 |\ga| + 2 -2n} dx.$$
 Moreover, by (\ref{i:Gexp}) and (\ref{i:v_exp})
  \begin{align}
  & \int_{\de B_{\gd} \cap \mathbb R^n_+} (\sum_{i=1}^n(
     v_{\gep}^2 \de_k h_{ik}+ \frac{n}{n-2} \de_k v_{\gep}^2 h_{ik})
   -\frac{4(n-1)}{n-2} \sum_{i=1}^n
  \gep^{\frac{n-2}{2}} (v_{\gep} \de_i
G - G \de_i v_{\gep})) \frac{x_i}{|x|} d\gs \notag\\
 & \leq \; - \gep^{n-2} \mathcal{I} (p, \gd) + C  \sum_{i, k=1}^n \sum_{2 \leq |\ga|\leq d}  |h_{ik, \ga}| \gd^{ |\ga|- n +2}
 \gep^{n-2}+ C \gep^{n-1} \gd^{-n+1}. \notag
 \end{align} From these the assertion follows.
 \epf
 We are ready to prove Theorem~\ref{t:main}.
 \bpf[Proof of Theorem~\ref{t:main}]
  Since $p \notin \mathcal{Z},$ we have $\sum_{i, k=1}^n \sum_{2 \leq |\ga|\leq d} |h_{ik, \ga}| > 0.$
  Thus, by Proposition~\ref{p:main},
  $$ \int_M (\frac{4(n-1)}{n-2} |d \phi_{(\gep, \gd)}|_g^2 + R_g \phi_{(\gep, \gd)}^2) dV_g \notag\\
    <
   \mathcal{\mathcal{Q}} (B, \de B) (\int_{\de M}
   \phi_{(\gep, \gd)}^{\frac{2(n-1)}{n-2}}d\gs_g)^{\frac{n-2}{n-1}} $$
   for $\gep > 0$ sufficiently small. This completes the proof.
 \epf

 Now we consider the case that $p \in \mathcal{Z}.$ We recall a
 result about $\mathcal{I} (p, \gd).$
 \begin{prop} \cite{Bchen09} Let $p \in \de M.$
 Suppose $p \in \mathcal{Z}.$

 (i) The limit $\lim_{\gd \rightarrow 0} \mathcal{I} (p, \gd)$
 exists.

 (ii)The doubling of $(M \setminus \{p\}, G^{\frac{4}{n-2}} g)$ has
 a well-defined mass which equals $\lim_{\gd \rightarrow 0} \mathcal{I} (p, \gd)$
 up to a positive factor.
 \end{prop}

 \bpf[Proof of Theorem~\ref{t:main2}]
Since $p \in \mathcal{Z},$ we have $\sum_{i, k=1}^n \sum_{2 \leq
|\ga|\leq d} |h_{ik, \ga}| = 0.$
  By Proposition~\ref{p:main},
    \begin{align}
    \int_M (\frac{4(n-1)}{n-2} |d \phi_{(\gep, \gd)}|_g^2 + R_g \phi_{(\gep, \gd)}^2) dV_g
    &\leq
   \mathcal{\mathcal{Q}} (B, \de B) (\int_{\de M}
   \phi_{(\gep,   \gd)}^{\frac{2(n-1)}{n-2}}d\gs_g)^{\frac{n-2}{n-1}}\notag \\
    & - \;\gep^{n-2} \mathcal{I} (p, \gd)
    + C  \gd^{2d+4-n} \gep^{n-2}+ C \gd^{-n+1} \gep^{n-1} \notag
   \end{align}
   for $0 < 2\gep \leq \gd.$
   By assumption $\lim_{\gd \rightarrow 0} \mathcal{I} (p, \gd) > 0,$ we may
   choose $\gd$ sufficiently small such that $\mathcal{I} (p, \gd) - C  \gd^{2d+4-n} > 0.$
   We next choose $0 < \gep < \frac{\gd}{2}$ sufficiently small
   such that $\mathcal{I} (p, \gd) - C  \gd^{2d+4-n} - C \gd^{-n+1} \gep>
   0.$  Then
  $$ \int_M (\frac{4(n-1)}{n-2} |d \phi_{(\gep, \gd)}|_g^2 + R_g \phi_{(\gep, \gd)}^2) dV_g \notag\\
    <
   \mathcal{\mathcal{Q}} (B, \de B) (\int_{\de M}
   \phi_{(\gep, \gd)}^{\frac{2(n-1)}{n-2}}d\gs_g)^{\frac{n-2}{n-1}}.$$
 \epf
 \vskip 1em
\noindent \textbf{\Large Appendix: An elliptic system in $\mathbb
R^n_+$}
 \vskip 1em

In the appendix, we solve  a boundary value problem for an elliptic
system in $\mathbb{R}^n_+.$

 Let $B_{\frac{1}{2}}$ be the ball of radius $\frac{1}{2}$ equipped with the flat
 metric $g.$ We denote by $\mathcal{X}$ the space of vector fields
 $V \in H^1 (B_{\frac{1}{2}})$ such that $\langle V, \nu \rangle = 0$ on $\de B_{\frac{1}{2}},$
 where $\nu$ is the unit outer normal on $\de B_{\frac{1}{2}}.$ We also denote by
 $\mathcal{Y}$ the space of trace-free symmetric  two-tensors on $B_{\frac{1}{2}}$
 of class $L^2.$ Let  $\mathcal{D}: \mathcal{X} \rightarrow \mathcal{Y}$
 be the conformal killing operator, which satisfies
  $$(\mathcal{D} V)_{ik} = V_{i, k} + V_{k, i} - \frac{2}{n} div V
  g_{ik}.$$

 By stereographic projection, $B_{\frac{1}{2}}$ is conformal to the
 hemisphere $\mathbb{S}^n_+$ with standard metric $g_c.$ The metric $g_c$ 
 satisfies $g_c = u^{\frac{4}{n-2}} g,$ where $u= (\frac{2}{1+ 4 |x|^2})^{\frac{n-2}{2}}$
 for $|x| \leq \frac{1}{2}.$ We may define similarly
   $\mathcal{X}^*$ the space of vector fields
 $V \in H^1 (\mathbb{S}^n_+)$ such that $\langle V, \nu \rangle = 0$ on
 $\de \mathbb{S}^n_+,$
 where $\nu$ is the unit outer normal on $\de \mathbb{S}^n_+$,
 $\mathcal{Y}^*$ the space of trace-free symmetric  two-tensors on $\mathbb{S}^n_+$
 of class $L^2$  and  $\mathcal{D}^*: \mathcal{X}^* \rightarrow \mathcal{Y}^*$
  the conformal killing operator on the hemisphere.
  Then it follows that $V \in H^1 (\mathbb{S}^n_+)$ if and only if
  $V \in H^1 (B_{\frac{1}{2}}),$ and $\mathcal{D}^* V = 0$ if and only
  if $\mathcal{D} V = 0.$
 \begin{lemma} $\ker \mathcal{D}$ is finite dimensional.
 \end{lemma}
  \bpf In \cite{Bchen09}, it was shown (after Lemma 21) that $\ker \mathcal{D}^*$ is finite
 dimensional. Then the assertion follows easily.
  \epf

 We now define
 $\mathcal{X}_0 = \{V \in \mathcal{X}: \langle V, U\rangle_{L^2 (B_{\frac{1}{2}})} = 0 \;\text{for
 all}\;
  U \in \ker \mathcal{D}\}.$

 \begin{lemma} \label{l:H^1} For all $V \in \mathcal{X}_0,$ it holds
  $\|V\|_{H^1 (B_{\frac{1}{2}})}^2 \leq C \|\mathcal{D} V\|_{L^2(B_{\frac{1}{2}})}^2,$ where $C= C(n).$
 \end{lemma}
 \bpf
 Suppose the inequality does not hold, then there exist a sequence
 of vector fields $V^{(j)} \in \mathcal{X}_0$ such that $\|V^{(j)}\|_{H^1(B_{\frac{1}{2}})} =
 1$ for all $j$ and $\|\mathcal{D} V^{(j)}\|_{L^2 (B_{\frac{1}{2}})} \rightarrow 0$ as $j \rightarrow \infty.$
 By passing to a subsequence, $V^{(j)}\rightharpoonup V^{(0)}$ weakly in $H^1 (B_{\frac{1}{2}})$ for some $V^{(0)} \in \mathcal{X}_0.$
 It follows that $\mathcal{D}V^{(0)} = 0,$ and as a result  $V^{(0)} = 0.$
 Notice that $V^{(j)}\rightarrow V^{(0)}$ strongly in $L^2 (B_{\frac{1}{2}}).$
 Thus, $\|V^{(j)}\|_{L^2 (B_{\frac{1}{2}})} \rightarrow 0.$ 
 Therefore,
  $\|V^{(j)}\|_{L^2 (\mathbb{S}^n_+)} \rightarrow 0.$   By
 \cite{Bchen09} Lemma 21, $\|V^{(j)}\|_{H^1 (\mathbb{S}^n_+)} \rightarrow 0$ as  $j \rightarrow \infty.$
 Hence, $\|V^{(j)}\|_{H^1 (B_{\frac{1}{2}})} \rightarrow 0$ as  $j \rightarrow \infty.$ This gives a
 contradiction.
 \epf

 \begin{prop} \label{p:system1}
 Let  $h$ be a two-tensor in $\mathcal{Y}.$
 Then there exists a unique vector field $V \in \mathcal{X}_0$ such that
 $\langle h - \mathcal{D} V, \mathcal{D} U \rangle_{L^2 (B_{\frac{1}{2}})} =0$ for
 all $U \in \mathcal{X}.$

  Moreover, $\|V\|_{H^1 (B_{\frac{1}{2}})}^2 \leq C \|h\|_{L^2 (B_{\frac{1}{2}})}^2,$
   where $C = C(n).$
 \end{prop}
  \bpf
  It follows by the same argument in \cite{Bchen09} Proposition 23,
  and Lemma~\ref{l:H^1} above that the minimizer of $\|h -\mathcal{D}V \|_{L^2 (B_{\frac{1}{2}})}^2$
  exists in $\mathcal{X}_0,$ which satisfies the required properties.
    \epf
  We now consider another conformal map. The ball $B_{\frac{1}{2}}$ is conformal to
   $\mathbb{R}^n_+ \cup \{\infty\}.$
   The metric $g$ satisfies $g= v^{\frac{4}{n-2}} \gd,$ where
  $$v = \left(\frac{1}{(1+ x_n)^2 + \sum_{1 \leq a \leq n-1} x_a^2)}\right)^{
  \frac{n-2}{2}}.$$

 \begin{prop} \label{p:system} Let  $h$ be a smooth trace-free symmetric two-tensor on $\mathbb{R}^n_+$ with compact support.
  Then there exists a smooth vector field $V$ on $\mathbb{R}^n_+$
  such that
       $$ \left\{  \begin{array}{ll}
  \sum_{k =1}^n \de_k [v^{\frac{2n}{n-2}} (h_{ik} - \de_i V_k - \de_k V_i + \frac{2}{n} div V \gd_{ik})]= 0 &  in \, \mathbb{R}^n_+  \\
  \de_n V_a- h_{an}=  0 & on \, \de \mathbb{R}^n_+\\
  V_n = 0 & on \, \de \mathbb{R}^n_+
  \end{array}\right.
    $$ for $i = 1, \cdots, n$ and $a= 1, \cdots, n-1.$
  Moreover,
    $$\int_{\mathbb{R}^n_+} v^{\frac{2(n+2)}{n-2}} |V|^2 dx \leq C \int_{\mathbb{R}^n_+} v^{\frac{2n}{n-2}} |h|^2 dx,$$
    where $C = C(n).$
 \end{prop}

\bpf By Proposition~\ref{p:system1}, there exists a smooth vector
field $V$ such that
   $$\int_{\mathbb{R}^n_+}  v^{\frac{2(n+2)}{n-2}} (h_{ik} - \de_i V_k - \de_k V_i
   + \frac{2}{n} div V \gd_{ik}) \de_k U_i dx= 0$$ for all $U \in \mathcal{X}$
    and $V_n = 0$ on $\de \mathbb{R}^n_+.$
   By elliptic regularity (\cite{Hormander85} pp.245-249), $V$ is
   smooth. Hence, $\sum_{k =1}^n \de_k [v^{\frac{2n}{n-2}} (h_{ik} - \de_i V_k - \de_k V_i + \frac{2}{n} div V \gd_{ik})]= 0 $
    on $\mathbb{R}^n_+$ and $\de_n V_a- h_{an}=0$ on $\de \mathbb{R}^n_+.$
\epf

 \begin{prop}
     Let  $h_{ik} = \eta (\frac{|x|}{\rho}) \sum_{|\ga|= 2}^d h_{ik, \ga}
     x^{\ga}$
     be a trace-free symmetric two-tensor,  where $d = [\frac{n-2}{2}]$,  $\rho \geq 1$
      and $\eta (t)$ be a fixed cut-off function which satisfies $\eta (t) =
      0$ for $t \geq 2.$
      Suppose $V$ is the vector field constructed in Proposition~\ref{p:system}.
      Then for  $x \in \mathbb{R}^n_+,$
    $$
    |\de^{\gb} V|^2 (x) \leq C(n, |\gb|) \sum_{i,k} \sum_{2 \leq |\ga| \leq d} |h_{ik, \ga}|^2 (1+ |x|)^{2|\ga|+2 - 2|\gb|}
    $$ for every multi-index $\gb.$
 \end{prop}
\bpf
   The proof is similar to \cite{Brendle07} Proposition 23 and Corollary 24.

      Without loss of generality we may assume
      $h_{ik} = \eta (\frac{|x|}{\rho}) \sum_{|\ga|= l} h_{ik, \ga} x^{\ga},$ where $2 \leq l \leq d.$
   We first prove that
    \begin{align} \label{i:system-reg}
    \sup_{r \geq 1} r^{-2l -n -2} \int_{(B_{2r} \setminus B_r) \cap \mathbb{R}^n_+} |V|^2 dx &\leq
    C\int_{\mathbb{R}^n_+} ((1+x_n)^2 + \sum_{a= 1}^{n-1} x_a^2)^{-n-2} |V|^2 dx \notag\\
     &+ C\sup_{r \geq 1} r^{-2l -n} \int_{(B_{2r} \setminus B_r) \cap \mathbb{R}^n_+} |h|^2 dx.
    \end{align}

   Suppose (\ref{i:system-reg}) does not hold, there exist  sequences
$h_{ik}^{(s)}$ and  $V^{(s)}$ such that
      $$\sup_{r \geq 1} r^{-2l -n -2} \int_{(B_{2r} \setminus B_r) \cap \mathbb{R}^n_+} |V^{(s)}|^2 dx = 1,$$
       $$\lim_{s \rightarrow \infty}\int_{\mathbb{R}^n_+} ((1+x_n)^2 + \sum_{a= 1}^{n-1} x_a^2)^{-(n+2)} |V^{(s)}|^2 dx =0,$$ and
       $$\lim_{s \rightarrow \infty}\sup_{r \geq 1} r^{-2l -n}  \int_{(B_{2r} \setminus B_r) \cap \mathbb{R}^n_+} |h^{(s)}|^2 dx = 0.$$
       Therefore, there exists a sequence $\rho^{(s)} \rightarrow \infty$ such that
       $$(\rho^{(s)})^{-2l -n-2} \int_{(B_{2 \rho^{(s)}} \setminus B_{\rho^{(s)}}) \cap \mathbb{R}^n_+}|V^{(s)}|^2 dx \geq \frac{1}{2}.$$
       Let
         $\tilde{h}_{ik}^{(s)} = (\rho^{(s)})^{- l} |x|^{-4} (|x|^2 \gd_{ij} - 2 x_i x_j) (|x|^2 \gd_{kl} -2 x_k x_l) h_{jl}^{(s)} (\frac{\rho^{(s)}
         x}{|x|^2}),$
        and $\tilde{V}_j = (\rho^{(s)})^{- l-1} (|x|^2 \gd_{ij} - 2 x_i x_j) V_i^{(s)} (\frac{\rho^{(s)}
        x}{|x|^2}).$
        Then they satisfy
         $$ \sum_{k =1}^n \de_k [((1 + \frac{x_n}{\rho^{(s)}})^2 + \sum_{a= 1}^{n-1} (\frac{x_a}{\rho^{(s)}})^2)^{-n} (\tilde{h}_{ik} - \de_i
         \tilde{V}_k-
          \de_k \tilde{V}_i + \frac{2}{n} div \tilde{V} \gd_{ik})]= 0$$ in $\mathbb
          R^n_+$ for $i = 1, \cdots, n,$ and
          $\tilde{V}_n =\de_n \tilde{V}_a -\tilde{h}_{an}=0$
          on $\de \mathbb{R}^n_+.$ Thus, by passing to a subsequence, $\tilde{V}_j$ converges  weakly to a vector field
        $V \in W^{1, 2}_{loc} (\mathbb{R}^n_+ \setminus \{0\}).$ $V$ satisfies
         $$\sum_{k =1}^n \de_k [ - \de_i V_k - \de_k V_i + \frac{2}{n} div V \gd_{ik}]= 0$$
         weakly in $\mathbb{R}^n_+ \setminus \{0\}$ for $i = 1, \cdots,
         n,$ and
         $V_n = \de_n V_a= 0$ on $\de \mathbb{R}^n_+ \setminus \{0\}.$
         By elliptic regularity theory, $V$ is smooth
         in $\mathbb{R}^n_+ \setminus \{0\}$.
        Thus, $V$  satisfies $\Delta V_j + \frac{n-2}{n} \de_j div V = 0.$ This implies $\Delta div V =0.$
        Moreover, on  $\de \mathbb{R}^n_+ \setminus \{0\}$ we have
        $0= \Delta V_n + \frac{n-2}{n} \de_n div V  = 2 \frac{n-1}{n} \de_n \de_n V_n.$
        Therefore,   $\de_n div V = 0$ on $\de \mathbb{R}^n_+ \setminus \{0\}.$
        We now define the function $div V$ on $\mathbb{R}^n \setminus
        \{0\}$ by  standard reflection. Then
         $div V$ is a $C^{2, 1}$ harmonic function in $\mathbb R^n \setminus \{0\}.$         Since $\sup_{\mathbb{R}^n \setminus \{0\}} |x|^{2l} |div V|^2$ is bounded, we obtain
         $div V = 0$ in $\mathbb{R}^n \setminus \{0\}.$  Thus, $\Delta V_j = 0$
        in $\mathbb{R}^n_+ \setminus \{0\}.$
        By the same reflection argument applied to the function $V_a,$
        we get $V_a$ is a $C^{2, 1}$ harmonic function in $\mathbb R^n \setminus \{0\}.$     Since  $\sup_{\mathbb{R}^n \setminus \{0\}} |x|^{2l-2} |V|^2 < \infty,$
          we have $V_a = 0$ in $\mathbb R^n \setminus \{0\}.$  Finally, since $\de_n V_n = div V= 0,$ using the same reflection argument again we obtain $V_n = 0$ in $\mathbb R^n \setminus \{0\}.$
          This contradicts to $\int_{(B_1 \setminus B_{\frac{1}{2}}) \cap \mathbb{R}^n_+} |V|^2 dx >
          0.$ Thus, (\ref{i:system-reg}) holds.

    Now since we have
    \begin{align}
     \int_{\mathbb{R}^n_+} ((1+x_n)^2 + \sum_{a= 1}^{n-1}x_a^2)^{-n-2} |V|^2 dx
    &\leq C \int_{\mathbb{R}^n_+} ((1+x_n)^2 + \sum_{a= 1}^{n-1}x_a^2)^{-n} |h|^2 dx\notag \\
    & \leq C \sum_{i, k= 1}^n \sum_{|\ga| = l} |h_{ik, \ga}|^2,
    \notag
    \end{align}
    then by (\ref{i:system-reg})
     $\sup_{r \geq 1} r^{-2l -n -2} \int_{(B_{2r} \setminus B_r) \cap \mathbb{R}^n_+} |V|^2 dx \leq C \sum_{i, k= 1}^n \sum_{|\ga| = l} |h_{ik,
     \ga}|^2.$
     Finally, by elliptic regularity
      $|\de^{\gb} V|^2 (x) \leq C(n, |\gb|) \sum_{i,k= 1}^n \sum_{|\ga| =l} |h_{ik, \ga}|^2 (1+ |x|)^{2|\ga|+2 - 2|\gb|}.$
  \epf

Institute for Advanced Study, Princeton, NJ
 \par
 Email address: \textsf{sophie@math.ias.edu}

\end{document}